\documentclass[12pt]{article}
\usepackage{setspace}
\usepackage{color}
\usepackage{lipsum}
\usepackage{verbatim}
\usepackage{amsfonts,amsmath,amssymb,mathrsfs}
\usepackage{graphicx}
\usepackage{epstopdf}
\usepackage{natbib}
\usepackage{algorithmic}
\usepackage{amsopn}
\usepackage{mathtools}
\usepackage{cases}
\usepackage{enumerate}
\usepackage{graphicx,subfigure}
\newtheorem{theorem}{Theorem}[section]
\newtheorem{lemma}{Lemma}[section]
\newtheorem{definition}{Definition}[section]
\newtheorem{remark}{Remark}[section]

\newtheorem{proposition}{Proposition}[section]

\usepackage{geometry}
\newcommand{\blue}{\textcolor{blue}}
\renewcommand{\epsilon}{\varepsilon}
\newcommand{\esssup}{\mathop{\operatorname{ess\,sup}}}

\numberwithin {equation} {section} 
\newcommand{\be}{\begin{equation}}
\newcommand{\ee}{\end{equation}}

\begin{document}
\title{\LARGE \bf Optimal Comfortable Consumption under Epstein-Zin utility}
\author{ Dejian Tian\thanks{School of Mathematics, China University of Mining and Technology, China. Email: djtian@cumt.edu.cn} \qquad
Weidong Tian\thanks{Corresponding author. Department of Finance, 
University of North Carolina at Charlotte, USA.  Email: wtian1@charlotte.edu. The authors are grateful to  Tao Pang and Jianfeng Zhang for their insightful comments and helpful discussions.} 
\qquad
Jianjun Zhou \thanks{College of Science, Northwest A\&F University, China. Email: zhoujianjun@nwsuaf.edu.cn.} \qquad
Zimu Zhu\thanks{Fintech Thrust, Hong Kong University of Science and Technology (Guangzhou), China. Email: zimuzhu@hkust-gz.edu.cn}
}

%
\date{}
\maketitle
\thispagestyle {empty} 

\newpage
\renewcommand{\abstractname}{
{\LARGE Optimal Comfortable Consumption under Epstein-Zin utility}
\\[1in] Abstract}

\begin{abstract}

We solve the optimal portfolio choice problem under Epstein--Zin utility with a time-varying consumption constraint, where closed-form expressions for neither the primal nor the dual value function are available. We establish the dynamic programming principle for the value function and prove that it is a viscosity solution of the corresponding Hamilton--Jacobi--Bellman equation. We further establish the $C^2$ regularity of the value function and derive a verification theorem using stochastic perturbation techniques. Finally, we provide an explicit characterization of the constrained region. The proposed methodology extends naturally to other constrained portfolio choice problems under the Epstein--Zin utility.

\vspace{0.9cm}



\end{abstract}
\thispagestyle{empty}

\newpage

\setcounter{page}{1}

\section{Introduction}

This paper studies the optimal portfolio and consumption problem under a dynamic consumption constraint in a continuous-time, infinite-horizon setting. Specifically, the consumption rate is required either to remain above a fixed level or to grow at least as fast as a prescribed benchmark rate, such as the inflation rate, at all times. Such constraints arise naturally in the management of retirement funds, endowment funds, and other long-horizon investment problems, as they ensure that investors maintain sufficient resources to finance essential consumption needs throughout the planning horizon. Consequently, they have been studied extensively under standard time-separable preferences; see, for example, \cite{BP92}, \cite{D95}, \cite{R91}, and \cite{PTT26}. This paper investigates the corresponding optimization problem under Epstein--Zin preferences and develops a novel stochastic control framework for its solution.

The Epstein–Zin preference allows for a separation between risk aversion and the elasticity of intertemporal substitution, making them particularly attractive for asset pricing, long-term portfolio choice, and models of economic growth; see, for instance, \cite{EZ89}, \cite{BY04}, \cite{GP15}, and \cite{B20}. However, its recursive structure, together with the non-Lipschitz and potentially non-convex or non-concave nature of the Epstein--Zin aggregator, creates analytical challenges in continuous time that are absent under expected utility or standard recursive utility models with Lipschitz aggregators, such as those studied in \cite{DE92} and \cite{EPQ97}. The presence of dynamic consumption constraints adds a further layer of complexity to the associated stochastic control problem.

In a recent paper,  \cite{HHJ23b} solve the corresponding optimal portfolio and consumption problem in the same market setting considered here, but without the consumption constraint, and therefore provide a natural unconstrained benchmark. In their setting, a closed-form expression for the value function is available, allowing the verification theorem to be established directly. In the constrained problem considered here, neither the value function nor the optimal policy admits such an explicit characterization. To the best of our knowledge, this constrained portfolio optimization problem under Epstein--Zin preferences has not previously been studied. We develop a stochastic control approach to analyze the problem, with five main contributions.

First, we establish the dynamic programming principle (DPP), \textit{Theorem \ref{theoremddp111}}, for the portfolio optimization problem under Epstein–Zin preferences. Since the DPP is closely related to BSDE representations of recursive preferences (see \cite{peng}), we develop new existence and uniqueness results for the associated BSDEs with infinite or random time horizon to handle the Epstein–Zin aggregator. We also derive key a priori estimates for the consumption and wealth processes, which are essential for the analysis of the BSDEs. The proposed approach to proving the DPP accommodates general path-independent constraints on the control variables.

Second, we prove in \textit{Theorem \ref{prop:vis}} that the value function is a viscosity solution of the associated Hamilton--Jacobi--Bellman (HJB) equation. The proof builds on the dynamic programming principle established above. Since the wealth process is not uniformly continuous in time, standard arguments for proving the viscosity subsolution property are not directly applicable. To overcome this difficulty, we adapt the approach of \cite{E08}. The viscosity characterization is then used to establish the $C^1$ regularity of the value function.

Third, we establish the $C^2$ regularity of the value function  in \textit{Theorem \ref{thm:recursive}}. 
We first establish the $C^1$ regularity and derive additional structural properties needed for the subsequent verification analysis. Our proof of the $C^2$  regularity draws on the Epstein–Zin utility literature, the viscosity solution characterization and $C^1$ regularity theory for HJB equations (see \cite{FS06, PZ18}), as well as a method introduced in \cite{XY16} for the time-separable utility case. 


Fourth, We characterize the value function through a verification theorem (\textit{Theorem \ref{th:verification}}). Building on the properties of the value function established in the previous three contributions, we identify a class of candidate value functions. The verification theorem then shows that this class contains a unique element, which coincides with the value function. The proof combines comparison arguments for stochastic differential equations with the stochastic perturbation techniques developed in \cite{HHJ21}, \cite{HHJ23b}, and \cite{PTT26}. A key ingredient is the establishment of a strong transversality condition, which is essential for the verification argument.


Finally, we explicitly characterize the constrained and unconstrained regions in \textit{Theorem \ref{th:two-region}}. Epstein–Zin preferences prevent the direct application of the classical approach in \cite{VZ97} and \cite{XY16}. Instead, we analyze the HJB equation separately in the two regions, express the optimal control in terms of the value function and its derivative, and exploit the structure of the resulting nonlinear ordinary differential equations.

The remainder of this paper is organized as follows. Section \ref{sec2} formulates the optimization problem under Epstein–Zin preferences. Section \ref{sec:DPP} establishes the dynamic programming principle (DPP), and Section \ref{sec:C2} establishes the $C^2$ regularity of the value function. Section \ref{sec:verification} presents the verification theorem and characterizes the constrained region. Finally, Section \ref{sec7} concludes the paper. All technical arguments are deferred to the appendices.

\subsection{Relevant Literature}

The methodology developed in this paper follows the viscosity solution approach for stochastic control problems in which the value function is unavailable in closed form. The key idea is to first establish, via the dynamic programming principle (DPP), that the value function is a viscosity solution of the associated HJB equation, and then derive its  $C^2$-regularity.  Early applications of this methodology include \cite{Z94} and \cite{E08}. More recently, \cite{FGG15} adopted a similar approach to establish the regularity of the value function for a class of non-standard preferences over both consumption and terminal wealth. 

Classical proofs of the DPP, such as those in \cite{yong22} and \cite{peng}, are not directly applicable in our setting because they require a careful analysis of backward stochastic differential equations (BSDEs) on infinite or random time horizons. We address this challenge by developing new technical arguments that accommodate unbounded control variables, state constraints, and the recursive structure of Epstein--Zin preferences. This yields a DPP for continuous-time portfolio optimization under Epstein--Zin preferences. Our approach differs substantially from the recent discrete-time analyses of \cite{BLV24} and \cite{JN24}, which establish the DPP by adapting the fixed-point arguments of \cite{BS20} and \cite{MM10} rather than through the analysis of the associated BSDEs.

For standard time-separable preferences, constrained portfolio choice has been extensively studied, and its distinguishing features relative to the unconstrained setting are well understood; see, for example, \cite{JK24}, \cite{SS15}, and \cite{Z94}. In contrast, under Epstein--Zin preferences, the complexity of the associated infinite-horizon BSDEs has limited most existing studies to unconstrained portfolio optimization, primarily in finite-horizon settings; see \cite{KSS13}, \cite{KSS17}, \cite{MX18}, \cite{MMS20}, \cite{SS03}, \cite{SS16}, \cite{X17}, \cite{hu24}, and \cite{FTZ26}. More recently, \cite{AH21} and \cite{AH23} extended this line of research to unconstrained problems in both finite- and infinite-horizon models.


A common feature of many existing portfolio optimization problems is the availability of an explicit value function, which allows the verification theorem to be proved directly; see, for example, \cite{D95}, \cite{SS99}, \cite{ET08}, \cite{KSS13}, \cite{DL11}, \cite{YK18}, \cite{MMS20}, and \cite{HHT26}. In the absence of an explicit representation, however, verification becomes significantly more difficult. Our principal methodological contribution is to develop a dynamic programming framework that establishes the verification theorem without requiring a closed-form expression for the value function.

Our work is also related to PDE-based approaches, particularly Perron's method, which provides an alternative route to verification without explicit solutions. For example, \cite{D21} studies an infinite-horizon Epstein--Zin portfolio optimization problem using Perron's method under restrictive assumptions on the market coefficients; see also \cite{JS12}, \cite{BZ15}, \cite{GH23}, and \cite{BL26}. By contrast, our approach is based on dynamic programming principle. It establishes the DPP under weaker assumptions, characterizes the value function through the associated HJB equation, and derives the regularity needed for verification.

\section{The optimization problem} \label{sec2}

 The probability space is defined as $\big(\Omega,(\mathcal{F}_t)_{t\geq0}, \mathcal{F}, P\big)$, where the information flow $(\mathcal{F}^t_s)_{s\geq t}$ is generated by a one-dimensional standard Brownian motion $(B_s-B_t)_{s\geq t}$, $\mathcal{F}_s = \mathcal{F}^0_s, \forall s  \ge 0$ and   $\mathcal{F}$ is the sigma-algebra generated by $(\mathcal{F}_t)_{t\geq0}$. Let ${\cal P}^t$ be the set of $(\mathcal{F}^t_s)_{s\geq t}$ progressively measurable process and ${\cal P}^t_{+}, {\cal P}^t_{++}$ the restriction of ${\cal P}^t$ to processes that take nonnegative and strictly positive values, respectively. When $t=0$, we drop the superscript $t$ and write these sets as ${\cal P}, {\cal P}_{+}, {\cal P}_{++}$.
 
  We consider a continuous-time financial market with one risky asset representing a stock index, whose price process evolves according to the dynamics
\begin{equation*}
d S_t=S_t(\mu d t+\sigma d B_t),
\end{equation*}
where $\mu > 0$ and $\sigma > 0$ are constants. Additionally, there is one risk-free asset with a constant rate of return $r$, and $0 < r < \mu$.

Let $\mathcal{A}(x)$ denote the set of \textit{feasible} consumption-investment strategies $(\pi_t, c_t)_{t\geq0}$ such that $\pi \in {\cal P}, c \in {\cal P}_{+}$,  the consumption rate $c_t \geq a, \forall t$, and there is a unique solution $(X_t)_{t\geq0}$ of the following stochastic differential equation 
\begin{equation}
\label{eq:wealth}
d X_t=r X_t d t+\pi_t X_{t}(\mu-r) d t+\pi_t X_{t} \sigma d B_t-c_t d t,
\end{equation}
with $X_0 = x$ and the non-negativity constraint $X_t \ge 0$ for all $t \ge 0$. Here, the investment strategy $\pi$ denotes the fraction of wealth  invested in the risky asset, and $X$ describes the wealth process. The positive number $a$ serves as a lower bound on the consumption rate, ensuring a minimum level of consumption throughout the investment horizon. In Section \ref{sec:verification} below, we will explain how the consumption constraint can be generalized to $c_t \ge a e^{\beta t}, \forall t \ge 0$.

Given the consumption constraint $c_{t}\geq a$, it is straightforward to see that the initial wealth must satisfy $x\geq \frac{a}{r}$.  Otherwise, there is no feasible consumption-investment strategy, implying that $\mathcal{A}(x)$ is an empty set for $ x< \frac{a}{r}$. Furthermore, if $x=\frac{a}{r}$, then the consumption rate must necessarily be $c_t = a $ for all time $t$. See \cite{ET08}.  In this case, the only feasible strategy is to invest entirely in the risk-free asset, as this ensures that the wealth remains sufficient to sustain the minimal consumption rate. 

Let $R$ and $S$ both lie in $(0,1) \cup (1,+\infty)$, and set $\mathbb{V}=(1-R)\mathbb{\bar{R}}_{+}$.  Let  
\begin{align*}
\nu=\frac{1-R}{1-S}, \quad  \rho=\frac{S-R}{1-R}=1-\frac{1}{\nu}.
\end{align*}
We define an aggregator as follows
\be
\label{eq:aggregator}
f(c,v)=\frac{c^{1-S}}{1-S}\big((1-R)v\big)^{\rho}, ~~c > 0, (1-R)v > 0, 
\ee
where $R$ denotes relative risk aversion and $S$ the agent's elasticity of intertemporal complementarity (EIC), defined as the reciprocal of the elasticity of intertemporal substitution. Let $\delta > 0$ denote the subjective discount rate, reflecting the agent's rate of time preference. 

The stochastic differential utility process $V^{c}$ associates to the consumption process $c \in {\cal P}_{+}$ and the aggregator $f$ is one that satisfies
\begin{equation}\label{eq:ez-utility}
V_{t}^{c}=\mathbb{E}\bigg[\int_{t}^{+\infty}e^{-\delta s} f(c_{s},V_{s}^{c})ds~\big|~\mathcal{F}_{t}\bigg],~~t\geq 0. 
\end{equation}
The aggregator (\ref{eq:aggregator}) differs slightly from the classical (minus) version in the Epstein-Zin literature; however, the corresponding differential utility processes retain all essential properties, as demonstrated in \cite{HHJ23a}. Moreover, the space $\mathscr{P}_{+}$ of consumption streams under the aggregator (\ref{eq:aggregator}) is broader than that of the classical aggregator, eliminating restrictive integral conditions (see, e.g., \cite{SS99}) that limit its applicability.

 In addition, \cite{HHJ23a, HHJ23b} demonstrate that the coefficients $R$ and $S$ must lie on the same side of unity (that is, $\nu > 0$) to ensure a well-defined utility process and prevent bubble formation. Furthermore, under the assumption $\nu\in(0,1)$,  Theorem 6.9 in \cite{HHJ23b} shows that for any consumption $c\in\mathscr{P}_{+}$, there exists a unique stochastic recursive utility process $V^{c}= (V_{t}^{c})_{t\geq0}$ that solves \eqref{eq:ez-utility}. 
    Given these rigorous considerations, we adopt the description of the Epstein-Zin utility formulation presented in \cite{HHJ23a}, and assume that  $\nu\in(0,1)$ (equivalently $0<S<R<1$, or $1<R<S$). It is worth noting, however, that \cite{HHJ25} and \cite{S26} study Epstein--Zin preferences under the alternative parameter regimes $\nu>1$ and $\nu<0$, respectively, each adopting a different formulation with distinct mathematical properties and economic interpretations.

The optimal portfolio choice problem is as follows
\begin{align}\label{eq:ez-problem}
J(x)=\sup _{(\pi, c) \in \mathcal{A}(x)} V_{0}^{c}=\sup _{c \in \mathcal{C}(x)} V_{0}^{c},
\end{align}where $V^{c}= (V_{t}^{c})_{t\geq0}$ denotes an appropriate solution of \eqref{eq:ez-utility}, and the consumption stream $c\in \mathcal{C}(x)$ if there is an investment process $\pi$ such that $(\pi,c)\in \mathcal{A}(x)$. We refer $J(x)$ as the (optimal) value function throughout the paper.


We consider $a=0$ as a benchmark model, as investigated by \cite{HHJ23a, HHJ23b}, which represents a case without constraints on consumption. The corresponding problem is as follows: 
\begin{align*}
J^{ez}(x)=\sup _{(\pi, c) \in \mathcal{A}_{ez}(x)} V_{0}^{c},~~~~~x>0,
\end{align*}where $\mathcal{A}_{ez}(x)$ corresponds to the feasible strategies set when  $a=0$. Theorem 8.1 in \cite{HHJ23b} establishes that 
the optimal strategies are given by
\begin{align*}
\pi^{ez}=\frac{\mu-r}{R\sigma^{2}}  \text{~~~~~and~~~~~}  c^{ez}=\eta X,\end{align*} and the optimal value function $J^{ez}(x)=\eta^{-\nu S}\frac{x^{1-R}}{1-R}$, under assumption that
\begin{align*}
\eta \equiv \frac{1}{S}\big[\delta+(S-1)\big(r+\frac{\kappa}{R}\big)\big] > 0,
\end{align*} 
where $\kappa=\frac{(\mu-r)^{2}}{2\sigma^{2}}$.

As another benchmark, where $R = S$ and thus $\nu = 1$, \cite{PTT26} solve the optimal portfolio choice problem in the time-separable case for any $a > 0$. It is still possible to derive an explicit value function under the condition $\eta > 0$. Therefore, throughout this paper, we assume that $\eta>0$ to ensure the existence of an optimal solution. Unlike the benchmark models, however, there is no closed-form expression for  $J(x)$ when $a> 0$, as will be shown below. As a result, the verification theorem cannot be established through the standard direct approach based on an explicit representation of the value function.

\section{Dynamic programming principle}
\label{sec:DPP}

In this section, we prove a Dynamic Programming Principle (DPP) under Epstein-Zin preferences in Theorem \ref{theoremddp111}. This DPP is a crucial step for proving the viscosity solution property of the value function in the next section.

The classical approach to establishing the DPP (see, e.g., \cite{yong22}) is not directly applicable in our setting for two key reasons. First, under recursive preferences, the problem is typically reformulated via an associated BSDE (see \cite{peng}).
However, in the Epstein--Zin framework with non-Lipschitz aggregator, the existence and uniqueness of solutions to the corresponding BSDE are not well studied yet. Second, standard proofs in the DPP literature —such as \cite{yong22}—depend heavily on unconstrained state spaces, whereas our setting involves both state and control constraints that complicate the argument.

To establish Theorem \ref{theoremddp111}, we first prove several technical lemmas for the corresponding SDEs in the Epstein-Zin framework. For every given bounded stopping time $\tau\geq0$ and random variable $\xi\geq0$, let $\mathcal{A}(\tau,\xi)$ denote the set of \textit{feasible} consumption-investment strategies $(\pi_s, c_s)_{s\geq \tau}$ such that $\pi \in {\cal P}^{\tau}$, $c \in {\cal P}^{\tau}_{+}$ and the consumption rate $c_s \geq a$, $\forall s\geq \tau$, and there is a unique solution $(X_s)_{s\geq \tau}$ of the following stochastic differential equation 
\begin{equation}
\label{eq:wealthmove}
d X_s=r X_s d s+\pi_s X_{s}(\mu-r) d s+\pi_s X_{s} \sigma d B_s-c_s d s,\ s\geq \tau, 
\end{equation}
with initial wealth $X_\tau = \xi$ and the non-negativity constraint $X_s \ge 0$ for all $s \ge \tau$. We also denote $X^{t,x,\pi,c}$ by the unique solution of \eqref{eq:wealthmove} with $(\tau,\xi)=(t,x)$.



\begin{lemma}\label{vtx}
Let $x\geq\frac{a}{r}$ and $t\geq0$.  For every $(\pi,c)\in \mathcal{A}(t,x)$, then there exists a unique process $(V_{s}^{t,x,\pi,c})$ satisfying
\begin{equation}\label{eq:ez-utility1}
V_{s}^{t,x,\pi,c}=\mathbb{E}\bigg[\int_{s}^{+\infty}e^{-\delta (l-t)} f(c_{l},V_{l}^{t,x,\pi,c})dl~\big|~\mathcal{F}^t_{s}\bigg],~~ s\geq t.
\end{equation}
Moreover,
\begin{align}\label{EVbounded}
\delta^{-\nu}\frac{a^{1-R}}{1-R}\leq V^{t,x,\pi,c}_t\leq \eta^{-\nu S}\frac{x^{1-R}}{1-R}.
\end{align}
\end{lemma}

We formulate the optimal portfolio choice problem conditional at time $t$  as follows.
\begin{align}\label{eq:ez-problemjia}
J(t,x)=\mathop{\esssup}\limits_{(\pi, c) \in \mathcal{A}(t,x)} V_{t}^{t,x,\pi,c}=\mathop{\esssup}\limits_{c \in \mathcal{C}(t,x)} V_{t}^{t,x,\pi,c},
~~~~x\geq\frac{a}{r},\end{align}
 where $V^{t,x,\pi,c}= (V_{s}^{t,x,\pi,c})_{s\geq t}$ denotes an appropriate solution of \eqref{eq:ez-utility1}, and the consumption stream $c\in \mathcal{C}(t,x)$ if there is an investment process $\pi$ such that $(\pi,c)\in \mathcal{A}(t,x)$. Clearly, $J(t,x)$ also satisfies \eqref{EVbounded}. For $t>0$ and random variable $\xi$ that is associated with $\mathcal{F}_t$, we define function $J(t,\xi):=J(t,x)|_{x=\xi}$.

Our approach to studying the process  $V^{t,x,\pi,c}$ is to characterize it as the solution to an infinite-horizon BSDE and then utilize the theory of infinite-horizon BSDEs. For this purpose, we consider the following infinite-horizon BSDE, for $ t\leq s\leq T<+\infty$:
\begin{align}\label{bsdetau}
Y^{t,x,\pi,c}_s=Y^{t,x,\pi,c}_T+\int^{T}_{s}f(c_l,Y^{t,x,\pi,c}_l)-\delta\nu Y^{t,x,\pi,c}_ldl-\int^{T}_{s}Z^{t,x,\pi,c}_ldB_l.
\end{align}
The existence and uniqueness of the solution to the above BSDE are established in the following lemma.

\begin{lemma}\label{th:finte-to-infinte}
Let $x\geq\frac{a}{r}$ and $t\geq0$. For every $(\pi,c)\in \mathcal{A}(t,x)$, 
 BSDE \eqref{bsdetau} has a unique solution $(Y,Z)$ with $Y^{t,x,\pi,c}_s=e^{\delta\nu (s-t)}V_s^{t,x,\pi,c}$,  where  $V_s^{t,x,\pi,c}$ is given in Lemma \ref{vtx}, and for every fixed $s\geq t$
\begin{align}\label{2602011}
\mathbb{E}[e^{-\delta\nu  (T-t)}Y_T^{t,x,\pi,c}|\mathcal{F}^t_s]
\rightarrow 0 \textrm{~~as~~} T\rightarrow +\infty.
\end{align}
Moreover, $\int^{+\infty}_0e^{-2\delta\nu (s-t)}|Z_s^{t,x,\pi,c}|^2ds< +\infty$.
\end{lemma}

Lemma~\ref{th:finte-to-infinte} essentially establishes a link between the infinite-horizon Epstein–Zin utility process and an infinite-horizon BSDE. \cite{DE92} obtained similar results under a Lipschitz assumption on the BSDE generator in \eqref{bsdetau}.  When $\nu<0$, \cite{D21} and \cite{S26} also studied the BSDE \eqref{bsdetau} for a class of consumptions in the Epstein-Zin framework. 

From the definition of $J(t,x)$ in \eqref{eq:ez-problemjia}, the value functional $J(0,x)$ is a deterministic function. Notice that $J(0,x) = J(x)$ defined in equation \eqref{eq:ez-problem}. The next result states the time-invariance property of $J(t,x)$, which will be used to formulate the DPP below. 

\begin{lemma}\label{25121201}
For every $(t,x)\in [0,+\infty)\times [\frac{a}{r},+\infty)$, we have $J(x)=J(t,x)$.
\end{lemma}

Since the DPP is related to the value function at any stopping time, 
we need to study the well-posedness of the finite-horizon BSDE for a given bounded stopping times $\tau\leq T$:
\begin{align}\label{eq:bsde-fix-terminal}
 Y_{t\wedge\tau}=\zeta+\int_{t\wedge\tau}^{\tau} f(c_{s},Y_{s})-\delta\nu Y_sds
 -\int^\tau_{t\wedge\tau}Z_sdB_s,\ \ 0\leq t\leq T.
\end{align}

Taken together, the following two lemmas establish the existence and uniqueness of a solution for this class of BSDEs.

\begin{lemma}\label{lem:bsde-rs<1}
 For any fixed finite horizon $T>0$ and bounded stopping times $\tau\leq T$, let $0<(1-R)\zeta\in \mathcal{F}_{\tau}$ and let $c \in {\cal L}_{+}$  satisfy $$\mathbb{E}\bigg[\big((1-R)\zeta\big)^{\frac{1}{\nu}}+\int_{0}^{\tau}c_{s}^{1-S}ds\bigg]<+\infty. $$ Then the finite-horizon BSDE
  \eqref{eq:bsde-fix-terminal}
admits a unique solution $(Y, Z)$. Moreover, $Y$ is strictly positive (resp., negative) when $0<R<1$ (resp., $R>1$), and of {class (D)}, and satisfy $\int_{0}^{\tau}Z_{s}^{2}ds<+\infty$. Furthermore, let   $(\zeta_{i}, c^{i})$, $i=1,2$, satisfies $\zeta_{1}\geq \zeta_{2}$ and $c^{1}\geq c^{2}$. Then $Y_{t}^{{1}}\geq Y_{t}^{{2}}$. In particular, if $c^{1}= c^{2}$, then $Y_{t\wedge\tau}^{{1}}- Y_{t\wedge\tau}^{{2}}\leq \mathbb{E}[\zeta_{1}-\zeta_{2}|\mathcal{F}_{t\wedge\tau}]$.
\end{lemma}

\begin{lemma}\label{Xbound}
Let $x\geq\frac{a}{r}$ and $0\leq t\leq T$. For each $(\pi,c)\in \mathcal{A}(t,x)$ and for any bounded stopping time $\tau\in[t,T]$ and $\alpha < 1$ with $\alpha \ne 0$, we have 
\begin{align*}
\mathbb{E}\bigg[X^{\alpha}_\tau+\int^\tau_tc^{\alpha}_sds\bigg]< +\infty.
\end{align*}
Moreover, $\mathbb{E}\big[e^{\int^\tau_t(\frac{\alpha(1-\alpha)}{4}\pi^2_s\sigma^2+\alpha X_s^{-1}c_s)ds}\big]<+\infty$.
\end{lemma}

\begin{remark}
 Even in the time-separable preference setting, many studies on portfolio choice impose certain integrability conditions on strategies to ensure the verification arguments rigorously (see the discussions in \cite{HHJ21} and \cite{PTT26}). Therefore, Lemma \ref{Xbound} is of independent interest, as it provides \emph{a priori} estimates for any admissible strategy  $(\pi,c)\in\mathcal{A}(t,x)$. 
\end{remark}



The final step in formulating the DPP is to define a family of backward semigroups, following the approach of \cite{peng}.
Given the initial value $x$,  a bounded stopping time $\tau>0$ and a control $(\pi,c)\in {\mathcal{A}}(x)$, we define
\begin{align}\label{gdpp}
G^{0,x,\pi,c}_{s,\tau}[J(\tau,X_\tau^{0,x,\pi,c})]:=\tilde{Y}^{0,x,\pi,c}_s,\ \ s\in[0,\tau],
\end{align}
where $(\tilde{Y}^{0,x,\pi,c}_s,\tilde{Z}^{0,x,\pi,c}_s)_{0\leq s\leq\tau}$ is the solution of the following
BSDE
\begin{equation}\label{bsdegdpp}
\begin{cases}
d\tilde{Y}^{0,x,\pi,c}_s =[\delta \nu \tilde{Y}^{0,x,\pi,c}_s-f(c_{s},\tilde{Y}^{0,x,\pi,c}_s)]ds +\tilde{Z}^{0,x,\pi,c}_sdB_s, \quad s\in [0,\tau], \\
\ \tilde{Y}^{0,x,\pi,c}_\tau=J(\tau,X_\tau^{0,x,\pi,c}).
\end{cases}
\end{equation}

The following lemma verifies that the definition in \eqref{gdpp} is well-defined. 

\begin{lemma}\label{finteBSDE}
For  every  bounded stopping time $\tau\in[0,T]$ with some constant $T>0$, and for every $(\pi,c)\in \mathcal{A}(x)$, the BSDE \eqref{bsdegdpp} admits a unique solution in the sense of Lemma \ref{lem:bsde-rs<1}. 
\end{lemma}

We now state the dynamic programming principle in our setting.

\begin{theorem}
\label{theoremddp111}(DPP) 
For any bounded stopping times $\tau>0$ and $x\geq \frac{a}{r}$,  the value function
$J$ satisfies the following  
\begin{align}\label{ddpG25}
J(x)=\mathop{\esssup}\limits_{(\pi,c)\in{\mathcal{A}}(x)}G^{0,x,\pi,c}_{0,\tau}[J(X^{0,x,\pi,c}_{\tau})].
\end{align}
 Alternatively, the DPP can be formulated as follows.
\begin{equation*}
J(x)
=
\operatorname*{ess\,sup}_{(\pi,c)\in\mathcal{A}(x)}
\mathbb{E}
\left[
\left.
\int_0^\tau e^{-\delta\nu s} f(c_s,Y_s)\,ds
+
e^{-\delta\nu\tau}
J\!\left(X^{0,x,\pi,c}_\tau\right)\,\right.
\right]
\end{equation*}
where $Y$ denotes the first component of the solution $(Y,Z)$ to the BSDE \eqref{bsdegdpp}, and $\tau$ is an arbitrary finite stopping time.
\end{theorem}

\begin{remark}
In our setting, admissible controls $(\pi,c)\in\mathcal{A}(x)$ must maintain non-negative wealth and satisfy the consumption constraint. These constraints prevent constructing the approximating controls $(\pi^{\varepsilon},c^{\varepsilon})\in\mathcal{A}(x)$ typically used in standard proofs of the DPP. To overcome this, our proof of the DPP follows the approach of \cite{fab1}, which relies on regular conditional probabilities to establish Lemma \ref{25121201} and Theorem \ref{theoremddp111}. This technique enables us to avoid relying on standard approximation arguments.
\end{remark}
\section{The $C^2$-regularity of the value function}\label{sec:C2} 

The main result of this section is Theorem \ref{thm:recursive}, which establishes 
the $C^2$-regularity of the value function. Our proof proceeds in three main steps. First, building on the DPP established in the previous section, we show that the value function is a viscosity solution. Next, we leverage this viscosity solution property to prove that the value function is  $C^1$-smooth. Finally, through a careful local analysis of the value function, we establish its $C^2$-regularity.

Specifically, we consider the viscosity solutions of 
\begin{equation}
\label{eq:Case2-vis}
-\delta\nu J(x)+\mathbf{H}\big(x,J(x),J'(x),J{''}(x)\big)=0,  \ x\in [\frac{a}{r},+\infty),
\end{equation}
where, for all $(x,k,p,q)\in
[\frac{a}{r},+\infty)\times \mathbb{R}_+\times \mathbb{R}\times \mathbb{R}$, 
\begin{align*}
\mathbf{H}(x,k,p,q)=\sup_{\alpha}[\alpha(\mu-r)p+\frac{1}{2}\sigma^2\alpha^2q]+\sup_{ c\geq a}[f(c,k)-cp]+rxp.
\end{align*}

We say $\varphi\in C_l^{2}\big((\frac{a}{r},+\infty)\big)$ if $\varphi\in C\big([\frac{a}{r},+\infty)\big)$ and $\varphi\in C^{2}([z,+\infty)) \text{~for all~} \ z>\frac{a}{r}$. For any $x\in [\frac{a}{r},+\infty)$ and $w \in C([\frac{a}{r},+\infty))$, define
\begin{align*}
  \mathcal{A}^+(x,w):=\Big\{\varphi\in C_l^{2}\big((\frac{a}{r},+\infty)\big)~:~ 
  0={w}(x)-\varphi(x)=\sup\nolimits_{y\in [a/r,+\infty)}({w}(y)- \varphi(y ))
\Big\}
\end{align*}
and
\begin{align*}
  \mathcal{A}^-(x,w):=\Big\{\varphi\in C_l^{2}\big((\frac{a}{r},+\infty)\big)~:~ 
  0={w}(x)+\varphi(x)=\inf\nolimits_{y\in [a/r,+\infty)}({w}(y)+\varphi(y))
\Big\}.
\end{align*}

\begin{definition}
 A function $w\in C\big([\frac{a}{r},+\infty)\big)$ is called a viscosity subsolution (resp.,  supersolution) to  \eqref{eq:Case2-vis} if whenever  $\varphi\in {\cal{A}}^+(x,w)$ (resp.,  $\varphi\in {\cal{A}}^-(x,w)$)  with $x>\frac{a}{r}$,  we have
\begin{align*}
 &-\delta\nu w(x)+{\mathbf{H}}\big(x,w(x),\varphi'(x),\varphi{''}(x)\big)\geq0,\\
\big(\text{resp.}, &-\delta\nu w(x)+{\mathbf{H}}\big(x,w(x),-\varphi'(x),-\varphi''(x)\big)\leq0\big).
\end{align*}
$w\in C\big([\frac{a}{r},+\infty)\big)$ is said to be a
viscosity solution to \eqref{eq:Case2-vis} if it is both a viscosity subsolution and a viscosity supersolution.
\end{definition}

\begin{theorem}\label{prop:vis}
The optimal value function $J$  is a continuous viscosity solution of \eqref{eq:Case2-vis}
in the class of concave functions with $J({a\over r})=\delta^{-\nu}{a^{1-R}\over 1-R}$.
\end{theorem}

The continuous and concavity property of the value function is proved in Proposition \ref{pro:ez-concave} in the Appendix. Typically, the viscosity solution property follows directly from the DPP in a standard stochastic control framework (see, e.g., \cite{FS06}). In our setting, however, due to the unboundedness of the control pair $(\pi,c)$, the wealth process $X^{0,x,\pi,c}$
 fails to be continuous in time uniformly over all admissible controls. Consequently, standard techniques are inapplicable. To overcome this difficulty, we adopt the approach of \cite{E08} and exploit the additional estimate to rigorously prove the viscosity subsolution property in Theorem \ref{prop:vis}.
      




The next key step is the $C^1$ smoothness of the value function.
Our proof is motivated by \cite{XY16} in the case of time-separable preferences; however, the Epstein--Zin setting requires establishing several technically demanding intermediate results (such as Propositions \ref{pro:ez-increase} and \ref{pro:ez-concave} in the Appendix) before the argument can be completed. Moreover, the proof relies crucially on the viscosity solution property of the value function established above.

\begin{proposition}
\label{prop:recursive-C1}
The optimal value function $J$ is $C^1(({a\over r},+\infty))$.
\end{proposition}

Building on the $C^1$-regularity  and viscosity solution property established above, the following result shows that $J$ is {\em strictly} increasing and  $J'$ is {\em strictly} decreasing.
We will subsequently use this result to perform the local analysis needed in the proof of Theorem \ref{thm:recursive}.

\begin{proposition}
\label{prop:strictly increasing}
(i) The value function $J$ is strictly increasing. (ii) For any point $x_0>{a\over r}$ where $J''(x_0)$ exists, we have $J''(x_0)<0$. Moreover, $J'$ is strictly decreasing on $(\frac{a}{r}, \infty)$.
\end{proposition}




By \eqref{eq:Case2-vis} and Proposition \ref{prop:recursive-C1}, the optimal consumption rate is given by
\begin{align*}
c^{*}=\max\Big\{a, (J')^{-1/S} \big((1-R)J\big)^{\frac{\rho}{S}}\Big\}.
\end{align*}
We define the {\em unconstrained domain} as
\begin{align}
\label{eq:U}
  \mathcal{U} =\Big\{x\in(\frac{a}{r},+\infty):  \big((1-R)J(x)\big)^{\frac{\rho}{S}} (J'(x))^{-\frac{1}{S}}>a\Big\},
\end{align}
and the {\em constrained domain} as 
\begin{align}
\label{eq:B}
 \mathcal{B} =\big\{x\in(\frac{a}{r},+\infty): \big((1-R)J(x)\big)^{\frac{\rho}{S}} (J'(x))^{-\frac{1}{S}}<a\Big\}.
\end{align}
By Proposition \ref{prop:recursive-C1}, both ${\mathcal U}$ and ${\mathcal B}$ are open sets. 


\begin{remark}
    In the time-separable preference case, where $R=S$ or $\rho = 0$, the constrained region is given by $\mathcal{B} = \left\{x: J'(x) > a^{-S} \right\}$. Because $J'$ is strictly decreasing, there exists a unique positive number $b > \frac{a}{r}$ such that $\mathcal{B} = (\frac{a}{r}, b)$, and thus  $\mathcal{U} = (b, +\infty)$ (see \cite{PTT26}). In contrast, for the Epstein-Zin utility, the structure of the constrained region depends on the joint behavior of $J$ and $J'$, making the characterization of both regions significantly more complex. We explore this intricate structure in detail in the next section.
\end{remark}


Finally, we establish the $C^2$-regularity of the value function.
By employing a standard transformation method together with Propositions \ref{prop:recursive-C1}--\ref{prop:strictly increasing}, we first show $C^2$-regularity in $\mathcal{U}$ and $\mathcal{B}$. We then complete the proof through a careful local analysis of the HJB equation near the boundary. 

\begin{theorem}
\label{thm:recursive}
The value function $J$ belongs to $C^2(({a\over r},+\infty))\cap C([{a\over r},+\infty))$. Moreover, the optimal investment strategy is given by
\begin{align*}
\pi^{*}=-\frac{\mu-r}{\sigma^{2}} \frac{J'(x)}{xJ''(x)}.
\end{align*}
\end{theorem}




\section{Verification theorem and the constrained region}
\label{sec:verification}

The main result of this section are verification theorem (Theorem \ref{th:verification}) and the explicit characterization of the constrained region (Theorem \ref{th:two-region}). We first introduce a class $\mathcal{J}$ of candidate value functions and show that the true value function belongs to this class. We then prove the verification theorem, which states that any candidate value function $J\in\mathcal{J}$ coincides with the value function. Consequently, $\mathcal{J}$ is a singleton.
We further show that the constrained region is a open interval $(\frac{a}{r}, \hat{x})$ where $\hat{x}$ is uniquely characterized by the verification theorem.

\begin{definition} A function $J:(\frac{a}{r}, +\infty) \rightarrow \mathbb{R} $ is said to be a {\em candidate value function} of the problem (\ref{eq:ez-problem}), if $J$ satisfies the following conditions:

\begin{itemize}
\item[(A1)]   $J$ is $C^{2}$ on $(\frac{a}{r}, +\infty)$; $J' > 0$ and $J''<0$; and $(1-R) J>0$.
\item[(A2)]   $J$ satisfies the following HJB equation 
\begin{align}\label{eq:hjb}
\delta\nu J(x)=\sup_{c\geq a, ~\alpha} \Big\{ f(c,J(x))+J'(x)[rx+\alpha(\mu-r)-c]+\frac{\alpha^{2}\sigma^{2}}{2}J''(x)\Big\}.
\end{align}
\item[(A3)]   $J$ satisfies the  boundary conditions at  $\frac{a}{r}$: 
\begin{equation*}
J(\frac{a}{r}):=\lim_{x \downarrow \frac{a}{r}}J(x) = \delta^{-\nu}{a^{1-R}\over 1-R} \text{~~and~~} J'(\frac{a}{r}) :=\lim_{x\downarrow \frac{a}{r}}J'(x)= +\infty.
\end{equation*}
\item[(A4)] For any $x > \frac{a}{r},$
$\delta^{-\nu} \frac{(rx)^{1-R}}{1-R} \le J(x) \le \eta^{-\nu S}\frac{x^{1-R}}{1-R}.$
\end{itemize}
\end{definition}


Combined with our previous results, the following proposition establishes that the value function is a candidate value function.

\begin{proposition}\label{pro:J-daoshu}
We consider the following two cases:
\begin{itemize}
    \item $R>1$ and $r-\delta+{S-1\over S}\kappa>0$,
    \item $0<R<1$.
\end{itemize}
Under these mild parameter assumptions,
the value function $J(\cdot)$ satisfies 
$J^{\prime}(\frac{a}{r})=\lim_{x\downarrow \frac{a}{r}}J^{\prime}(x)=+\infty.$ 
\end{proposition}

When the context is clear, we use $J$ to  denote either the value function or the candidate value function. For verification purposes, it is essential to characterize the set $\mathcal{J}$ and identify one candidate value function within $\mathcal{J}$ that is the optimal value function. The main technical challenge under Epstein--Zin preferences, however, stems from the highly non-linear nature of the HJB equation \eqref{eq:hjb}, which complicates establishing the uniqueness of the candidate value function.

We derive several technical properties of the candidate value function to prove the verification theorem. First, we show that the candidate value function is $C^3$-smooth and that the its derivative is of order $x^{-R}$.

\begin{proposition}\label{pro-jprime-guji}
Suppose $J\in\mathcal{J}$. Then $J\in C^{3}$ and there exist three (not depending on $x$) positive constants $k_1$ and $k_2$, and $x_0>\frac{a}{r}$ such that 
$J'(x)\geq k_1 x^{-R}$ for all $x>a/r$, and $J'(x)\leq k_2 x^{-R}$ for all $x>x_0$.
\end{proposition}

Second, we analyze the corresponding strategies generated by the candidate value function $J \in \mathcal{J}$.
For any candidate value function $J$, solving the HJB equation \eqref{eq:hjb}  yields the candidate optimal feedback strategy $(\pi,c)$, defined by 
\begin{align}\label{eq:h-strategy}
\pi(x):=-\frac{\mu-r}{\sigma^{2}} \frac{J^{\prime}(x)}{x J^{\prime\prime}(x)}I_{x>\frac{a}{r}}  \textrm{~~~and~~~}   c(x)=\max\Big\{a, (J'(x))^{-1/S} \big((1-R)J\big)^{\frac{\rho}{S}}I_{x>\frac{a}{r}}\Big\}.
\end{align}
The corresponding feedback wealth process $X$ satisfies the following SDE:
\begin{align*}
\left\{\begin{array}{ll}
d X_t=r X_t dt+  \Big(\pi (X_{t})X_t(\mu-r) -c(X_{t}) \Big)d t+\pi(X_{t})X_t\sigma d B_t, \\
X_{0}=x>\frac{a}{r}.
\end{array}
\right.
\end{align*}
 In the following proposition, we show that the wealth process corresponding to a candidate value function remains strictly above $a/r$. This guarantees well-behaved feedback strategies $(\pi, c)$ and the corresponding wealth process.

\begin{proposition}\label{pro: SDE}
Let $J$ be a candidate value function. For any $x>\frac{a}{r}$, then the corresponding feedback strategy $(\pi, c)\in \mathcal{A}(x)$ and wealth process $X_t>\frac{a}{r}$ for all $t\geq0$.  
\end{proposition}

We now present the verification theorem for our setting. Its proof relies on three key ingredients. First, we show that the value function is the {\em smallest} candidate value function by extending the stochastic perturbation method developed in \cite{HHJ21, HHJ23b, PTT26}. Second, we prove a {\em strong transversality condition} for the candidate value function. Finally, the result follows from standard verification arguments together with Propositions \ref{pro-jprime-guji}--\ref{pro: SDE}.

\begin{theorem}
\label{th:verification}
    Assume that $x>\frac{a}{r}$. Let $J$ be a candidate value function satisfying (A1)-(A4), and $(\pi^*, c^*)$ be the corresponding feedback strategy induced by $J$, and $X^*$ be the corresponding wealth process. Then, the candidate value function $J$ is 
the optimal value function for the problem \eqref{eq:ez-problem}, and the corresponding optimal strategy is given by $(\pi^*, c^*)$. That is,
\begin{align*}
J(x)=V_0^{c^*} =\mathbb{E}\bigg[\int_{0}^{+\infty}e^{-\delta s}f(c^*(X^*_{s}), V_s^{c^*})ds\bigg].
\end{align*}
\end{theorem}




We next explicitly characterize the constrained and unconstrained regions, as detailed in Theorem \ref{th:two-region}. 

\begin{theorem}\label{th:two-region}
We consider the following two cases:
\begin{itemize}
    \item $R>1$ and $r-\delta+{S-1\over S}\kappa>0$,
    \item $0<R<1$ and $\kappa+\rho\delta\nu>0$.
\end{itemize}
Under these assumptions, there exists a unique critical wealth level $\hat{x} > \frac{a}{r}$ such that $\mathcal{B}=({a\over r}, \hat{x})$ and $\mathcal{U}=(\hat{x},+\infty)$.
\end{theorem}

The condition in Theorem~\ref{th:two-region} is mild, essentially requiring the Sharpe ratio of the risky asset to be reasonably large. Specifically, when $\kappa > \frac{S}{S-1}(\delta - r)$ when $R > 1$, and $\kappa > \delta (1-\nu)$ when $0 < R < 1$. Under the same condition, the next result follows from Theorem \ref{th:verification} - Theorem \ref{th:two-region}, which states that the value function $J$ and $\hat{x}$ are uniquely determined by a system of ordinary differential equations.


\begin{proposition}
Assume $\kappa > \frac{S}{S-1}(\delta - r)$ when $R > 1$, and $\kappa > \delta (1-\nu)$ when $0 < R < 1$. Then, there exists a {\em unique} strictly increasing and strictly concave function $J(x)$ over $[\frac{a}{r}, \infty)$ and a {\em unique} number $\hat{x} > \frac{a}{r}$ satisfying: $(1-R)J(x) >0, J(\frac{a}{r}) = \delta^{-\nu}{a^{1-R}\over 1-R},  J'(\frac{a}{r}) = +\infty$, $\delta^{-\nu} \frac{(rx)^{1-R}}{1-R} \le J(x) \le \eta^{-\nu S}\frac{x^{1-R}}{1-R}, \forall x$, along with
\begin{equation*}
-\kappa {(J')^2\over J''}+{S\over 1-S}((1-R)J)^{\rho\over S}(J')^{1-{1\over S}}+rxJ'-\delta\nu J=0, ~~ \ \ x > \hat{x};
\end{equation*}
and  
\begin{equation*}
-\kappa {(J')^2\over J''}+f(a,J)-aJ'+rxJ'-\delta\nu J=0, ~\ \ \frac{a}{r} < x < \hat{x}.
\end{equation*}
Moreover, this function $J$ is the value function of the optimal portfolio choice problem \eqref{eq:ez-problem} and the constrained region for this problem is given by $(\frac{a}{r}, \hat{x})$.
\end{proposition}

To finish the discussion of this section, we demonstrate how our model can be adapted to accommodate a general consumption constraint of the form $c_t \ge a e^{\beta t}$, $\forall t\geq0$, where the real number $\beta$ represents the minimal growth rate of the consumption, for example, due to inflation. This specification allows for a growing consumption floor over time while preserving the tractability of the model.
 	
Let $\tilde{c}_t = c_t e^{-\beta t}$. With this transformation, the original utility function based on consumption $c_t$ can be rewritten as
\begin{equation*}
V_{t}^{\tilde{c}}=\mathbb{E}\bigg[\int_{t}^{+\infty}e^{-\tilde{\delta} s} f(\tilde{c}_{s},V_{s}^{\tilde{c}})ds~\big|~\mathcal{F}_{t}\bigg],~~t\geq 0,
\end{equation*}where 
     $\tilde{\delta} = \delta + \beta(S-1)$.
	The wealth equation (\ref{eq:wealth}) for $X_t$ is equivalent to the following process for $\tilde{X}_t = e^{-\beta t}X_t$ as follows
	\begin{equation*}
	d\tilde{X}_t = (r - \beta) \tilde{X}_t dt + \tilde{\alpha}_t  (\mu-r) dt - \tilde{c}_t dt + \tilde{\alpha}_t \sigma dB_t,~~~~ \tilde{\alpha}_t = \alpha_t e^{-\beta t}.   
	\end{equation*}
	Thus, we can reduce the portfolio choice problem under the constraint $c_t \ge a  e^{\beta t}$ into a portfolio choice problem (\ref{eq:ez-problem}) in which $r$ is replaced by $\tilde{r} = r - \beta$, and the subjective discount factor $\delta$ is replaced by $\tilde{\delta}$. Moreover, by setting $\tilde{\mu} - \tilde{r} = \mu - r$, then $\tilde{\mu} = \mu - \beta$. Therefore, we establish the following proposition.

\begin{proposition}
\label{prop:extension}Suppose $x>\frac{a}{r-\beta}>0$ and $\eta>0$. 
		Let $\tilde{J}(x)$ be the candidate value function in Theorem \ref{th:verification} with parameters $\{\tilde{\mu}, \tilde{r}, \tilde{\delta}\}$. Then $\tilde{J}(x)$ is the value function of the optimal portfolio choice problem under the constraint $c_t \ge a e^{\beta t}$. Moreover, the optimal consumption-investment strategy, $\tilde{c}^*_t = e^{\beta t} c^*_t$, $\tilde{\alpha}^*_t = e^{\beta t} \alpha^*_t$, $\forall t\geq0$, determined in $\tilde{J}(x)$  leads to the optimal consumption-investment strategy for the constraint $c_t \ge a e^{\beta t}$, $t\geq0$.
	\end{proposition}

\section{Conclusions}\label{sec7}


We solve an optimal portfolio choice problem under Epstein–Zin preferences and a time-varying minimum consumption constraint. We establish the dynamic programming principle and show that the value function is a classical $C^2$ solution to the associated HJB equation. By synthesizing stochastic analysis, martingale techniques, and stochastic perturbations, we prove a verification theorem that uniquely determines the value function without an a priori candidate. Lastly, we explicitly characterize the constrained and unconstrained regions.
  
Beyond the specific portfolio framework analyzed here, our approach establishes a versatile foundation for constrained stochastic control under Epstein–Zin utility. More broadly, it resolves technical challenges in continuous-time recursive utility models characterized by non-Lipschitz, nonconcave, or nonconvex aggregators. Future research may employ this framework to investigate portfolio optimization problems with other classes of constraints.


\newpage

\newpage

\renewcommand {\theequation}{A-\arabic{equation}} \setcounter
{equation}{0}
\renewcommand {\thelemma}{A.\arabic{lemma}} \setcounter
{theorem}{0}
\renewcommand {\theproposition}{A.\arabic{proposition}} \setcounter
{theorem}{0}
\renewcommand {\theremark}{A.\arabic{remark}} \setcounter
{theorem}{0}
\setcounter{equation}{0}
\section*{Appendix A: Proofs in Section \ref{sec:DPP}}
\label{sec:DPP-appendix}

We begin by establishing two useful properties of the value function.

\begin{proposition}\label{pro:ez-increase}  The value function $J$  is finite-valued,  increasing on  $[a/r, +\infty)$,  and $(1-R)J(x)>0$ for all $x\geq a/r$, and $ J\big(\frac{a}{r}\big) =\delta^{-\nu}\frac{a^{1-R}}{1-R}$. Moreover, for all $x \ge \frac{a}{r}$, the following inequality holds
\begin{align}
    \label{eq:bound}
    \delta^{-\nu} \frac{(rx)^{1-R}}{1-R} \le J(x) \le \eta^{-\nu S}\frac{x^{1-R}}{1-R}.
\end{align}
Besides, 
$\lim_{x\rightarrow +\infty}J(x)=+\infty$ when $R\in(0,1)$, and $\lim_{x\rightarrow +\infty}J(x)=0$ when $R>1$.
\end{proposition}
\noindent\textbf{Proof.}  
Since $J(x)\leq J^{ez}(x)$, $J$ is finite-valued. It is also evident that $J$ is increasing. When $x=\frac{a}{r}$, the set $\mathcal{C}(x)$ contains only one element $c_{t}=a$ for all $t\geq0$.  In this case, the stochastic differential utility satisfies the ordinary differential equation $$d\bar{V}^a_t=-\nu a^{1-S}e^{-\delta t}(\bar{V}^a_t)^{\rho}dt,~~~\bar{V}^a_\infty=0$$ with $\bar{V}^a_t=(1-R)V^a_{t}$. Solving this equation yields $V^a_t=\delta^{-\nu}e^{-\delta\nu t}\frac{a^{1-R}}{1-R}$. This result implies 
\begin{align}\label{eq:J-ar}
J\big(\frac{a}{r}\big)=V^a_{0}=\delta^{-\nu}\frac{a^{1-R}}{1-R}.\end{align}


We next consider a feasible strategy $(\pi, c) = (0, rx)$ and the corresponding stochastic differential utility process $V_t^{(rx)}$. Using the same approach as in the derivation of equation (\ref{eq:J-ar}), we obtain
\begin{align*}
V_t^{(rx)} = \delta^{-\nu} e^{-\delta \nu t} \frac{(rx)^{1-R}}{1-R},
\end{align*}
which serves as the lower bound for $J(x)$ in equation (\ref{eq:bound}) at $t=0$. It follows from the inequality \eqref{eq:bound} that $\lim_{x\rightarrow +\infty}J(x)=+\infty$ when $R\in(0,1)$, and $\lim_{x\rightarrow +\infty}J(x)=0$ when $R>1$. Finally, equation (\ref{eq:bound}) implies that $(1-R)J(x)>0$, $\forall x\geq \frac{a}{r}.$
\hfill$\Box$



\begin{proposition}\label{pro:ez-concave}  The optimal value function $J$ is  concave and uniformly continuous on the domain $[a/r, +\infty)$.  
\end{proposition}
\noindent\textbf{Proof.} 
For any $c^{1}\in \mathcal{C}(x_{1})$ and $c^{2}\in \mathcal{C}(x_{2})$ and $0 < \lambda < 1$, let $c^{\lambda} \equiv \lambda c^{1}+(1-\lambda)c^{2}$. We first establish that
$V_{t}^{c^{\lambda}}\geq \lambda V_{t}^{c^{1}}+(1-\lambda) V_{t}^{c^{2}}, \forall t\geq0. $
Note that by Proposition \ref{pro:ez-increase}, $V_0^c\geq V_0^a$, and  $(V^c_t)_{t\geq0}$ is uniformly integrable for any given $c\in\mathcal{C}(x)$.   
We proceed by analyzing the two cases separately.

{\em Case 1:  $S>1$ and $R>1$.}  In this case, $f(c, v)$ is joint concave in $(c,v)$. Then, for any $t\geq0$,  we have
\begin{align*} 
&V_{t}^{c^{\lambda}}-[\lambda V_{t}^{c^{1}}+(1-\lambda) V_{t}^{c^{2}}]\\
=&\mathbb{E}\left[\int_{t}^{\infty}e^{-\delta s}\left( f(c^{\lambda}_{s}, V_{s}^{c^{\lambda}})-\lambda  f(c^{1}_{s}, V_{s}^{c^{1}})- (1-\lambda) f(c^{2}_{s}, V_{s}^{c^{2}})\right)ds~\big|~\mathcal{F}_{t}\right]\\
\geq & \mathbb{E}\left[\int_{t}^{\infty}e^{-\delta s}\left( f(c^{\lambda}_{s}, V_{s}^{c^{\lambda}})- f(c^{\lambda}_{s}, \lambda V_{s}^{c^{1}}+(1-\lambda) V_{s}^{c^{2}}) \right)ds~\big|~\mathcal{F}_{t}\right]\\
\geq & \mathbb{E}\left[\int_{t}^{\infty} e^{-\delta s} \frac{\partial f}{\partial v}(c^{\lambda}_{s}, V_{s}^{c^{\lambda}}) \left( V_{s}^{c^{\lambda}}-[\lambda V_{s}^{c^{1}}+(1-\lambda) V_{s}^{c^{2}}] \right) ds~\big|~\mathcal{F}_{t}\right].
\end{align*}
Given that  $\frac{\partial f}{\partial v}(c^{\lambda}_{s}, V_{s}^{c^{\lambda}})\leq 0$, it follows from Lemma C3 (remains valid when $T=+\infty$) in \cite{SS99} 
that  $V_{t}^{c^{\lambda}}-[\lambda V_{t}^{c^{1}}+(1-\lambda) V_{t}^{c^{2}}]\geq 0$ for all $t\geq0$.  

{\em Case 2: $R$, $S\in (0,1)$.} In this situation, $f(c, v)$ is no longer joint concave in $(c,v)$. 
Motivated by Proposition 2.4 of \cite{X17},  we approach the proof by applying a transformation method as follows. For each  $i=1,2,\lambda\in(0,1)$, since $\mathbb{E}\left[\int_{0}^{\infty}e^{-\delta s} f(c^{i}_{s}, V_{s}^{c^{i}}) ds\right]<+\infty$, 
 it follows that $V^{c^{i}}_{t}$ is strictly positive for all $t \in (0,+\infty)$ and  $V^{c^{i}}_{\infty}=0$.  By the representation theorem (Chapter 3, Corollary 3, p189) of \cite{P05}, we have
$$V_{t}^{c^{i}}=\int_{t}^{+\infty}e^{-\delta s}f(c^{i}_{s}, V_{s}^{c^{i}})ds-\int_{t}^{+\infty}Z_{s}^{i}dB_{s},$$where $\int_{0}^{+\infty}(Z_{s}^{i})^{2}ds<+\infty$ for $i=1,2,\lambda \in(0,1)$. Define $\tilde{V}_{t}^{c^{i}} :=(V_{t}^{c^{i}})^{1-\rho}$ and $\tilde{Z}^{i}_{t}:=(1-\rho)(V_{t}^{c^{i}})^{-\rho}Z_{t}^{i}$. Applying  It\^{o}'s formula yields
$$\tilde{V}_{t}^{c^{i}}=\int_{t}^{+\infty}F(s, c^{i}_{s}, \tilde{V}_{s}^{c^{i}}, \tilde{Z}_{s}^{i})ds-\int_{t}^{+\infty}\tilde{Z}_{s}^{i}dB_{s},$$where
$F(t,c,\tilde{v}, \tilde{z})=\frac{(1-\rho)(1-R)^{\rho}}{1-S}e^{-\delta t}c^{1-S}+\frac{\rho}{2(1-\rho)}\frac{|\tilde{z}|^{2}}{\tilde{v}}$.  Since $\rho<0$ and $S, R\in(0,1)$, the function  $F(t,c,\tilde{v}, \tilde{z})$ is jointly concave in $(c,\tilde{v}, \tilde{z})$.   

Let $\Delta \tilde{V} :=\lambda \tilde{V}^{c^{1}}+(1-\lambda) \tilde{V}^{c^{2}}$ and $\Delta \tilde{Z} :=\lambda \tilde{Z}^{1}+(1-\lambda) \tilde{Z}^{2}$. Then, we have
$$\Delta \tilde{V}_{t}=\int_{t}^{+\infty} (F(s,c_{s}^{\lambda}, \Delta \tilde{V}_{s}, \Delta \tilde{Z}_{s})+A_{s})ds-\int_{t}^{+\infty}\Delta \tilde{Z}_{s}dB_{s},$$where $A_{s}=\lambda F(s, c^{1}_{s}, \tilde{V}_{s}^{c^{1}}, \tilde{Z}_{s}^{1})+(1-\lambda)F(s, c^{2}_{s}, \tilde{V}_{s}^{c^{2}}, \tilde{Z}_{s}^{2})-F(s,c_{s}^{\lambda}, \Delta \tilde{V}_{s}, \Delta \tilde{Z}_{s})\leq0$.  Moreover, applying It\^o formula gives
$$(\Delta \tilde{V}_{t})^{\frac{1}{1-\rho}}=\mathbb{E}\left[\int_{t}^{+\infty}e^{-\delta s}f(c_{s}^{\lambda}, (\Delta \tilde{V}_{s})^{\frac{1}{1-\rho}})+\frac{(\Delta \tilde{V}_{s})^{\frac{\rho}{1-\rho}}}{1-\rho}A_{s}ds~\big|~\mathcal{F}_{t}\right],$$
Since $A_s \le 0$, we obtain
\begin{align*}
V_{t}^{c^{\lambda}}-(\Delta \tilde{V}_{t})^{\frac{1}{1-\rho}}\geq &\mathbb{E}\left[\int_{t}^{+\infty}e^{-\delta s}\left( f(c^{\lambda}_{s}, V_{s}^{c^{\lambda}})-f(c_{s}^{\lambda}, (\Delta \tilde{V}_{s})^{\frac{1}{1-\rho}})\right)ds~\big|~\mathcal{F}_{t}\right]\\
\geq &\mathbb{E}\left[\int_{t}^{+\infty}e^{-\delta s} \frac{\partial f}{\partial v}(c_{s}^{\lambda}, (\Delta \tilde{V}_{s})^{\frac{1}{1-\rho}}) \left(V_{s}^{c^{\lambda}}-(\Delta \tilde{V}_{s})^{\frac{1}{1-\rho}}\right)ds~\big|~\mathcal{F}_{t}\right],
\end{align*}
where the second inequality follows from the fact that $f(c,v)$ is convex in $v$.  Moreover, since $\frac{\partial f}{\partial v}\leq0$, applying Lemma C3 in \cite{SS99} together with $\rho<0$ yields
$$V_{t}^{c^{\lambda}}\geq(\Delta \tilde{V}_{t})^{\frac{1}{1-\rho}}\geq \lambda (\tilde{V}_{t}^{c^{1}})^{\frac{1}{1-\rho}}+(1-\lambda) (\tilde{V}_{t}^{c^{2}})^{\frac{1}{1-\rho}}=\lambda V_{t}^{c^{1}}+(1-\lambda) V_{t}^{c^{2}},$$
by virtue of the inequality $(\lambda x+(1-\lambda) y)^{p}\geq \lambda x^{p}+(1-\lambda) y^{p}$ for any $p,\lambda\in(0,1)$ and $x,y\geq0$.
Taking the supremum over all feasible consumption strategies gives $J(\lambda x_{1}+ (1-\lambda)x_{2})\geq \lambda J( x_{1})+ (1-\lambda)J(x_{2})$, 
thereby establishing the concavity of $J$.

The continuity of $J$ on $({a\over r},\infty)$ follows from the concavity property established above. We next prove that $\lim_{x\downarrow \frac{a}{r}} J(x)=J({a\over r})$ by contradiction. If not, then there exist $x_n\downarrow {a\over r}$ and $\epsilon_0>0$ such that $J(x_n)\geq J({a\over r})+\epsilon_0,\forall n$. Furthermore, there exist strategies $(c^n,\pi^n)\in\mathcal{A}(x_n)$ such that 
\begin{equation}
\label{eq:continuity-contradiction}
V_0^{c^n}>J({a\over r})+{\epsilon_0\over 2},~~\forall n.
\end{equation}
 From the wealth process, we have
\begin{align*}
(X_t^n-{a\over r})e^{-rt}+\int_0^{t}e^{-rs}(c_s^n-a)ds
=(x_n-{a\over r})+\int_0^{t}e^{-rs}\pi_s^nX_s^n\sigma d\tilde{B}_s,
\end{align*}where $d\tilde{B}_s=dB_s+\frac{\mu-r}{\sigma}ds$.
Since the left-hand side of the above equation is non-negative, an argument similar to (8.21)–(8.23) in Chapter 5 of \cite{KS91} shows that as $x_n\downarrow {a\over r}$,  $c_s^n\rightarrow a$, a.s. (taking a subsequence if necessary). On the other hand, let $V_0^a$ denote the utility at $0$ for $X_0={a\over r}$ and $c_t\equiv a$. We have
\begin{align*}
0\leq V^{c^n}_0-V^a_0&=\mathbb{E}\left[\int_0^{\infty}e^{-\delta s}(f(c_s^n,V_s^{c^n})-f(a,V_s^a))ds\right]\\
&\leq\mathbb{E}\left[\int_0^{\infty}e^{-\delta s}(f(c_s^n,V_s^{c^n})-f(a,V_s^{c^n}))ds\right]\\
&=\mathbb{E}\left[\int_0^{\infty}e^{-\delta s}\frac{(c_s^n)^{1-S}-a^{1-S}}{1-S} ((1-R)V_s^{c^n})^{\rho}
ds\right],
\end{align*}where the inequality follows from the facts that $f(c,v)$ is decreasing in $v$. For sufficient large $n$, we have $$((1-R)V_s^{c^n})^{\rho}\leq ((1-R)V_s^{a})^{\rho}+((1-R)V_s^{a+1})^{\rho},~~s\geq0,$$ where $V_s^a=\delta^{-\nu}\frac{a^{1-R}}{1-R}e^{-\delta\nu s}$ and $V_s^{a+1}=\delta^{-\nu}\frac{(a+1)^{1-R}}{1-R}e^{-\delta\nu s}$. Thus, the dominated convergence theorem implies that $V_0^{c^n}\rightarrow V_0^a=J({a\over r})$, which  contradicts (\ref{eq:continuity-contradiction}). Therefore, $J$ is continuous, and since it is concave, it is also uniformly continuous on $[\frac{a}{r}, +\infty)$. \hfill$\Box$



\subsection*{Proof of Lemma \ref{vtx}}

For any $t\geq0$, define $\hat{\mathcal{F}}_s:=\mathcal{F}^t_{s+t}$, ${B}^t_s:=B_{t+s}-B_t$, $s\geq0$. Let $\hat{\mathcal{A}}(x)$ denote the set of feasible consumption-investment strategies with $(\Omega, \mathcal{F}, (\hat{\mathcal{F}}_s)_{s\geq0}, \mathbb{P})$ and ${B}^t$. 
For every fixed $(\pi,c)\in \mathcal{A}(t,x)$, define
$(\hat{\pi}_s,\hat{c}_s):=(\pi_{s+t},c_{s+t}), ~~s\geq0.
$
It is clear that $(\hat{\pi},\hat{c})\in \hat{\mathcal{A}}(x)$. Define $\hat{V}_{s}^{x,\hat{\pi},\hat{c}}:=V_{s+t}^{t,x,\pi,c}$, $s\geq0$, equation \eqref{eq:ez-utility1} becomes
\begin{align}
\hat{V}_{s}^{x,\hat{\pi},\hat{c}}=V_{s+t}^{t,x,\pi,c}&=\mathbb{E}\left[\int_{s+t}^{\infty}e^{-\delta (l-t)} f(c_{l},V_{l}^{t,x,\pi,c})dl~\big|~\mathcal{F}^t_{s+t}\right]\nonumber\\
&=\mathbb{E}\left[\int_{s}^{\infty}e^{-\delta l} f(\hat{c}_{l},\hat{V}_{l}^{x,\hat{\pi},\hat{c}})dl~\big|~\hat{\mathcal{F}}_{s}\right].\label{eq:ez-utility10512}
\end{align}
By Theorem 6.9 in \cite{HHJ23b}, \eqref{eq:ez-utility10512} has a unique solution $\hat{V}^{x,\hat{\pi},\hat{c}}$. Consequently, \eqref{eq:ez-utility1} has a unique solution $V^{t,x,{\pi},{c}}$ on $[t,+\infty)$. Moreover, Proposition \ref{pro:ez-increase} yields that \eqref{EVbounded} holds. \hfill$\Box$

\subsection*{Proof of Lemma \ref{th:finte-to-infinte}}
By Lemma \ref{vtx}, for $
\mathbb{E}[|V^{t,x,\pi,c}_t|]<+\infty,$
for any $(\pi,c)\in\mathcal{A}(t,x)$. 
Then, by the martingale representation theorem, there exists a unique $(\mathcal{F}^t_s)_{s\geq t}$-adapted process  $(Z^1_s)_{s\geq t}$ such that $\int^{\infty}_{t}|Z^1_s|^2ds<+\infty$
and for all $s\geq t$, 
$$
\mathbb{E}\left[\int_{t}^{\infty}e^{-\delta (l-t)} f(c_{l},V_{l}^{t,x,\pi,c})dl~\big|~\mathcal{F}^t_{s}\right]
=\mathbb{E}\left[\int_{t}^{\infty}e^{-\delta (l-t)} f(c_{l},V_{l}^{t,x,\pi,c})dl\right]+\int^s_tZ^1_ldB_l.
$$In other words, 
\begin{align*}
&\mathbb{E}\left[\int_{s}^{\infty}e^{-\delta (l-t)} f(c_{l},V_{l}^{t,x,\pi,c})dl~\big|~\mathcal{F}^t_{s}\right]\\
=&\mathbb{E}\left[\int_{t}^{\infty}e^{-\delta (l-t))} f(c_{l},V_{l}^{t,x,\pi,c})dl\right]+\int^s_tZ^1_ldB_l
-\int_{t}^{s}e^{-\delta (l-t)} f(c_{l},V_{s}^{t,x,\pi,c})dl.
\end{align*}
Equivalently,
$$V_s^{t,x,\pi,c}-V_T^{t,x,\pi,c}=\int_{s}^{T}e^{-\delta (l-t)} f(c_{l},V_{l}^{t,x,\pi,c})dl
 -\int^T_sZ^1_ldB_l,\ \ t\leq s< T<+\infty.
$$
Applying It\^o's formula to $e^{\delta \nu (s-t)}V_s^{t,x,\pi,c}$,
one has $$de^{\delta \nu (s-t)}V_s^{t,x,\pi,c}=[\delta \nu e^{\delta \nu (s-t)}V_s^{t,x,\pi,c}-e^{-\delta (s-t)}e^{\delta \nu (s-t)}f(c_{s},V_{s}^{t,x,\pi,c})]ds
 +e^{\delta \nu (s-t)}Z^1_sdB_s.
$$
Then
$(Y_s^{t,x,\pi,c}, Z^{t,x,\pi,c}_s)_{s\geq t}:=(e^{\delta\nu  (s-t)}V_s^c, e^{\delta \nu (s-t)}Z^1_s)_{s\geq t}$
is the  solution of infinite-horizon BSDE \eqref{bsdetau}. Moreover, for every fixed $s\geq t$, taking $T\rightarrow +\infty$ yields
\begin{align*}
\mathbb{E}[(e^{-\delta\nu  (T-t)}Y_T^{t,x,\pi,c}|\mathcal{F}^t_s]
=\mathbb{E}[V_T^{t,x,\pi,c}|\mathcal{F}^t_s]
=\mathbb{E}\left[\int_{T}^{\infty}e^{-\delta (l-t)} f(c_{l},V_{l}^{t,x,\pi,c})dl~\big|~\mathcal{F}^t_{s}\right]\rightarrow 0.
\end{align*}
Thus, we obtain \eqref{2602011}.
 
 The uniqueness of $(Y^{t,x,\pi,c}, Z^{t,x,\pi,c})$ follows from the uniqueness of $V^{t,x,\pi,c}$. In fact, assume that $(Y^{1,c}, Z^{1,c})$ is another solution of BSDE (\eqref{bsdetau}). Applying It\^o's formula to $e^{-\delta \nu(s-t)}Y_s^{1,c}$ and then taking conditional expectation with respect to $\mathcal{F}^t_s$, for $t\leq s< T<+\infty$, one has
 $$e^{-\delta \nu (s-t)}Y^{1,c}_s-\mathbb{E}\left[e^{-\delta \nu (T-t)}Y^{1,c}_T|\mathcal{F}^t_s\right]=\mathbb{E}\left[\int_{s}^{T}e^{-\delta (l-t)} f(c_{l},e^{-\delta \nu (l-t) }Y^{1,c}_l)dl|\mathcal{F}^t_s\right].
 $$
 Letting $T\rightarrow +\infty$, from \eqref{2602011}, we have 
  $$
  e^{-\delta \nu (s-t)}Y^{1,c}_s=\mathbb{E}\left[\int_{s}^{\infty}e^{-\delta (l-t)} f(c_{l},e^{-\delta \nu (l-t) }Y^{1,c}_l)dl|\mathcal{F}^t_s\right],\ \ t\leq s<+\infty,
 $$
 and $Y^{1,c}_s=e^{\delta\nu  (s-t)}V_s^{t,x,\pi,c}=Y^{t,x,\pi,c}_s$ and  $Z^{1,c}_s=Z^{t,x,\pi,c}_s$.\hfill$\Box$

\subsection*{Proof of Lemma \ref{25121201}}
We denote by $\mathcal{P}^{B}_t$ the sigma field generated by the sets of the form $(s,r]\times A$, $t\leq s\leq k <+\infty$, $A\in \mathcal{B}^t_s$ and $\{t\}\times A$, $A\in \mathcal{B}^t_t$. Here $\mathcal{B}^t_s=\sigma(\mathcal{W}(l): t\leq l\leq s)$, $\Lambda_t=\{\omega \in C([t,+\infty);\mathbb{R}): \omega(t)=0\}$ and $
 \mathcal{W}: [t,+\infty)\times \Lambda\rightarrow \mathbb{R},\quad \mathcal{W}(s,\omega)=\omega(s).$

For each $\varepsilon>0$, there exists $(\pi,c)\in \mathcal{A}(x)$ such that
$J(x)\leq V^{x,\pi,c}_0+\varepsilon.$
By Lemma 2.20 in \cite{fab1},  there exists a $\mathcal{P}^{B}_0$-measurable function $g:[0,+\infty)\times \Lambda_0\rightarrow \mathbb{R}\times \mathbb{R}_+$ such that
$({\pi},{c})(s,\omega)=g(s, B(\cdot,\omega)),\ \mbox{for}\ \omega\in \Omega, \ s\geq0.$
Define
$$(\pi^1,c^1)(s,\omega):=g(s, B^t(\cdot,\omega)),\ \mbox{for}\ \omega\in \Omega, \ s\geq0,$$
and $$ (\pi^2,c^2)(s,\omega)= (\pi^1,c^1)(s-t,\omega)=g(s-t, B^t(\cdot,\omega)), \ \ s\geq t,$$
where
$B^t(s,\omega)=B(s+t,\omega)-B(t,\omega)$, $s\geq0.$
By its construction,$(\pi^2,c^2)\in \mathcal{A}(t,x)$, and 
$\mathcal{L}_{\mathbb{P}}({\pi},{c},B)= \mathcal{L}_{\mathbb{P}}(\pi^1,c^1,B^t).$
Moreover, because of the uniqueness solution of weak solution of BSDE (\ref{bsdetau}), $Y_0^{0,x,{\pi},{c}}=Y_t^{t,x,\pi^2,c^2}.$ 
Therefore,
\begin{align*}
J(x)\leq V^{x,\pi,c}_0+\varepsilon=Y^{0,x,{\pi},{c}}_0+\varepsilon=Y^{t,x,{\pi}^2,{c}^2}_t+\varepsilon=V^{t,x,\pi^2,c^2}_t+\varepsilon\leq J(t,x)+\varepsilon.
\end{align*}
Letting $\varepsilon\rightarrow0$, we have that $J(x)\leq J(t,x)$.

On the other hand, for every $\varepsilon>0$, 
 there exists $(\pi,c)\in \mathcal{A}(t,x)$ which is $({\mathcal{F}}^t_s)_{s\geq t}$ progressively measurable such that
$ J(t,x)\leq V^{t,x,\pi,c}_t+\varepsilon.$
By Lemma 2.20 in \cite{fab1}, there exists a $\mathcal{P}^{B}_t$-measurable function $g:[t,+\infty)\times \Lambda_t\rightarrow \mathbb{R}\times \mathbb{R}_+$ such that
$$ (\pi,c)(s,\omega)=g(s, B(\cdot,\omega)-B(t,\omega)),\ \mbox{for}\ \omega\in\Omega, \ s\geq t.$$
Define
$$(\pi^1,c^1)(s,\omega):=g(s, B(\cdot-t,\omega)),\ \mbox{for}\ \omega\in \Omega, \ s\geq t,$$
and
$$(\pi^2,c^2)(s,\omega):=(\pi^1,c^1)(s+t,\omega)=g(s+t, B(\cdot,\omega)),\ \mbox{for}\ \omega\in \Omega, \ s\geq 0.$$
 Notice that  $(\pi^2,c^2)\in \mathcal{A}(x),$ and $
\mathcal{L}_{\mathbb{P}}({\pi}^1,{c}^1,B)= \mathcal{L}_{\mathbb{P}}(\pi,c,B^t). $
Moreover, ${Y_0^{0,x,\pi^2,c^2}=Y_t^{t,x,\pi,c}.} $
 Therefore,
\begin{align*}
J(t,x)\leq V^{t,x,\pi,c}_t+\varepsilon=V^{x,\pi^2,c^2}_0+\varepsilon\leq J(x)+\varepsilon.
\end{align*}
Letting $\varepsilon\rightarrow0$, one has $J(t,x)\leq J(x)$. \hfill$\Box$

\subsection*{Proof of Lemma \ref{lem:bsde-rs<1}}

First, we consider the case $0<R<1$. 
Following the approach of Lemma 8.1 in \cite{hu24}, we construct the generator as follows:
 $$f_{m}(c,y)=\frac{(c\wedge m)^{1-S}}{1-S}\left((1-R)y\vee \frac{1}{m}\right)^{\rho},$$
 where $m > 0$. Under this construction,  $f_{m}$ is globally Lipschitz continuous in $y$. Moreover,
 $0\leq f_{m}\leq \frac{m^{1-S-\rho}}{1-S}$.
 Then the truncated BSDE
 \begin{align*}
 Y_{t\wedge\tau}^{c,m}=\zeta\wedge m+\int_{{t\wedge\tau}}^{\tau}f_{m}(c_{s},Y_{s}^{c,m})-\delta\nu Y^{c,m}_sds-\int^\tau_{t\wedge\tau}Z^{c,m}_sdB_s,\ \ 0\leq t\leq T,
\end{align*} admits a unique solution $(Y^{c,m}, Z^{c,m})$. Besides, $f_{m}$ and $\zeta\wedge m$ are both increasing with respect to $m$.  By the comparison theorem,  it follows that $Y^{c,m}$ is increasing with $m$, and nonnegative (taking generator and terminal random value as $0$).

We now establish an upper bound for $Y^{c,m}$ that is independent of $m$.
Applying It\^o's formula to $(Y^{c,m})^{\frac{1}{\nu}}$, taking the conditional expectation, we obtain
\begin{align*}
(Y_{t\wedge\tau}^{c,m})^{\frac{1}{\nu}}&\leq\mathbb{E}\left[ (\zeta\wedge m)^{\frac{1}{\nu}}+\int_{t\wedge\tau}^{\tau}\frac{1}{\nu} (Y_{s}^{c,m})^{-\rho} f_{m}(c_{s}, Y_{s}^{c,m})ds~\big|~\mathcal{F}_{t\wedge\tau}\right]\\
&\leq \mathbb{E}\left[ \zeta^{\frac{1}{\nu}}+(1-R)^{\rho-1}\int_{t\wedge\tau}^{\tau}c_{s}^{1-S}ds ~\big|~\mathcal{F}_{t\wedge\tau}\right].
\end{align*}
Hence,
$$Y_{t\wedge\tau}^{c,m}\leq  \left(\mathbb{E}\left[ \zeta^{\frac{1}{\nu}}+(1-R)^{\rho-1}\int_{t\wedge\tau}^{\tau}e^{-\delta s}c_{s}^{1-S}ds ~\big|~\mathcal{F}_{t\wedge\tau}\right]\right)^{\nu}=:\bar{Y}^{c}.$$
Then, by the localization argument in \cite{bh06}, we  obtain the existence of a solution to the BSDE \eqref{eq:bsde-fix-terminal}.
 Finally, since $f(c,y)$ is decreasing in $y$,uniqueness follows from an argument analogous to that of Proposition 2.2 in \cite{X17}.

Next, for the case $R>1$, we construct the generator  as follows: 
$$f_{m}(c,v)=\frac{(c\vee \frac{1}{m})^{1-S}}{1-S}\left(((1-R)v)\vee \frac{1}{m})\right)^{\rho}, \forall m > 0. $$
Then $f_{m}$ is globally Lipschitz continuous in $v$, and,
 $\frac{m^{S-1-\rho}}{1-S}\leq f_{m}\leq 0$.
 Then the truncated BSDE
 \begin{align*}
 V_t^{c,m}=\zeta\vee (-m)+\int_{t}^{T}e^{-\delta s} f_{m}(c_{s},V_{s}^{c,m})ds
 -\int^T_tZ^{c,m}_sdW_s,\ \ 0\leq t\leq T,
\end{align*}admits a unique solution $(V^{c,m}, Z^{c,m})$. Besides, $f_{m}$ and the terminal random value are both decreasing with respect to $m$,  by comparison theorem,  it implies $V^{c,m}$ is also decreasing with $m$, and nonpositive (taking generator and terminal random value as 0).

Similar to the case $0<R<1$, it suffices to establish a lower bound for $V^{c,m}$ that is independent of $m$.
Taking the conditional expectation, one has that
\begin{align*}
&((1-R)V_{t}^{c,m})^{1-\rho}
\leq\mathbb{E}\left[ ((1-R)(\zeta \vee (-m)))^{1-\rho}~\big|~\mathcal{F}_{t}\right]\\
+&\mathbb{E}\left[ \int_{t}^{T}  \nu(1-\rho)((1-R)V_{s}^{c,m})^{-\rho} e^{-\delta s}(c_{s}\vee \frac{1}{m})^{1-S} (((1-R)V_{s}^{c,m})\vee \frac{1}{m})^{\rho}ds
~\big|~\mathcal{F}_{t}\right]\\
\leq &\mathbb{E}\left[ ((1-R)\zeta)^{1-\rho}+\nu(1-\rho) \int_{t}^{T}e^{-\delta s}c_{s}^{1-S}ds ~\big|~\mathcal{F}_{t}\right].
\end{align*}
Thus, it implies that
\begin{align*}
V_{t}^{c,m}&\geq \frac{1}{1-R}\left(\mathbb{E}\left[ ((1-R)\zeta)^{1-\rho}+\nu(1-\rho) \int_{t}^{T}e^{-\delta s}c_{s}^{1-S}ds ~\big|~\mathcal{F}_{t}\right]\right)^{\nu}=:{\underline{V}}_{t}^{c}.
\end{align*}Similar to the case of $R<1$, the existence result can be obtained by the localization argument and the uniqueness result follows from the decreasing property of $f(c,y)$ in $y$. 
\hfill$\Box$

\subsection*{Proof of Lemma \ref{Xbound}}
The lemma is evident for $\alpha < 0$
 since the wealth process and consumption process are both bounded below in the current setting. We now assume that $\alpha\in(0,1)$ and consider stopping time $\tau\in[t,T]$. Direct calculation yields
\begin{align}\label{eq:x-tau}
X_\tau^{\alpha}=x^{\alpha}e^{\alpha\int^\tau_t(r+\pi_s(\mu-r)-\frac{1}{2}\pi_s^2\sigma^2-X_s^{-1}c_s)ds+\alpha\int^\tau_t\pi_s\sigma dB_s}.
\end{align}
By completing the square in the drift, the exponent in \eqref{eq:x-tau} can be decomposed as
\begin{align}\label{eq:x-zhishu-xiang}
&\alpha\int^\tau_t(r+\pi_s(\mu-r)-\frac{1}{2}\pi_s^2\sigma^2-X_s^{-1}c_s)ds+\alpha\int^\tau_t\pi_s\sigma dB_s\nonumber\\
=&\int^\tau_t-\frac{1}{2}(\alpha\pi_s\sigma)^2ds+\int^\tau_t\alpha\pi_s\sigma dB_s+\alpha\int^\tau_t\big(\frac{1}{2}(\mu-r)^2(1-\alpha)^{-1}\sigma^{-2}+r\big)ds\nonumber\\
&-\int_t^\tau \alpha X_s^{-1}c_sds-\frac{\alpha}{2}\int^\tau_t((1-\alpha)^{\frac{1}{2}}\pi_s\sigma+(\mu-r)(1-\alpha)^{-\frac{1}{2}}\sigma^{-1})^2ds.
\end{align}
Since the stochastic exponential term satisfies
\begin{align*}
\mathbb{E}\Big[\exp\left(\int^\tau_t-\frac{1}{2}(\alpha\pi_s\sigma)^2ds+\int^\tau_t\alpha\pi_s\sigma dB_s\right)\Big]\leq 1,
\end{align*}
dropping the last two non-positive terms in \eqref{eq:x-zhishu-xiang} yields $$\mathbb{E}[X_\tau^{\alpha}]\leq x^{\alpha}\exp\left(\alpha\int^T_t (r+\frac{1}{2}\frac{(\mu-r)^2}{(1-\alpha)\sigma^2})ds\right)=  x^{\alpha}\exp\left((r+\frac{\kappa}{1-\alpha})(T-t)\alpha\right).$$

On the other hand, for any $s\geq t$, let $$Y_s:=x+\int^s_tX_u[r+\pi_u(\mu-r)]du+\int^s_tX_u\pi_u\sigma dB_u,$$then $Y_s\geq X_s\geq \frac{a}{r}$. Using the same estimation procedure for $X_\tau$, we obtain 
\begin{align*}
\mathbb{E}[Y_\tau^{\alpha}]=&x^{\alpha}\mathbb{E}\exp\left(\alpha\int^\tau_t\big([r+\pi_s(\mu-r)]\frac{X_s}{Y_s}-\frac{1}{2}\pi_s^2\sigma^2\frac{X^2_s}{Y^2_s}\big)ds+\alpha\int^\tau_t\pi\sigma \frac{X_s}{Y_s}dB_s\right)\\
\leq& x^{\alpha}\exp\left((r+\frac{\kappa}{1-\alpha})(T-t)\alpha\right).
\end{align*}
Notice that for any $\tau\in[t,T]$, 
$\int^\tau_tc_sds=Y_\tau-X_\tau\leq Y_\tau.$ Thus, applying Hölder's inequality yields
\begin{align*}
\mathbb{E}\left[\int^\tau_tc^{\alpha}_sds\right]\leq \mathbb{E}\left[\int^T_tc^{\alpha}_sds\right]
\leq (T-t)^{(1-\alpha)}\mathbb{E}\left[\int^T_tc_sds\right]^{\alpha} \leq (T-t)^{(1-\alpha)}\mathbb{E}[Y_T^{\alpha}]<+\infty.
\end{align*}

Finally, returning to the equations \eqref{eq:x-tau} and \eqref{eq:x-zhishu-xiang} and moving the two negative terms to the left-hand side, we derive 
\begin{align*}
&X_\tau^{\alpha}\exp\bigg(\frac{\alpha}{2}\int^\tau_t[((1-\alpha)^{\frac{1}{2}}\pi_s\sigma+(\mu-r)(1-\alpha)^{-\frac{1}{2}}\sigma^{-1})^2+2X_s^{-1}c_s]ds\bigg)\\
=&x^{\alpha}\exp\bigg(\int^\tau_t-\frac{1}{2}(\alpha\pi_s\sigma)^2ds+\int^\tau_t\alpha\pi_s\sigma dB_s+\alpha\int^\tau_t\big(\frac{1}{2}(\mu-r)^2(1-\alpha)^{-1}\sigma^{-2}+r\big)ds\bigg).
\end{align*}
Using the lower bound $X_\tau\geq \frac{a}{r}$ and taking the expectations on both sides,  we obtain  
\begin{align*}
&\left(\frac{a}{r}\right)^{\alpha}\mathbb{E}\exp\left(\frac{\alpha}{2}\int^\tau_t [((1-\alpha)^{\frac{1}{2}}\pi_s\sigma+(\mu-r)(1-\alpha)^{-\frac{1}{2}}\sigma^{-1})^2+2X_s^{-1}c_s]ds\right)\\
\leq&x^{\alpha}\exp\left((r+\frac{\kappa}{1-\alpha})(T-t)\alpha\right).
\end{align*}
Consequently, 
\begin{align*}
\mathbb{E}\big[e^{\int^\tau_t\frac{\alpha(1-\alpha)}{4}\pi^2_s\sigma^2+\alpha X_s^{-1}c_sds}\big]\leq\left(\frac{rx}{a}\right)^{\alpha}e^{(r+\frac{\kappa}{1-\alpha})(T-t)\alpha}<+\infty.\end{align*}
The proof is complete. \hfill$\Box$

\subsection*{Proof of Lemma \ref{finteBSDE}}
 By Proposition \ref{pro:ez-increase} together Lemmas \ref{25121201} and \ref{Xbound}, for $0 < S < 1$,
 \begin{align*}
 \mathbb{E}[((1-R)J(T,X_{\tau}^{0,x,\pi,c}))^{\frac{1}{\nu}}]&=\mathbb{E}[((1-R)J(X_{\tau}^{0,x,\pi,c}))^{\frac{1}{\nu}}] \leq
 (\delta^{-1}r^{1-S}\vee\eta^{-S})\mathbb{E}[(X_{\tau}^{0,x,\pi,c})^{1-S}]<+\infty.
 \end{align*}
 On the other hand, by Lemma \ref{vtx}, for $S > 1$, the above inequality remains valid since  $X_{\tau}^{0,x,\pi,c}$ is bounded from below.
 Moreover, by Lemma \ref{Xbound} (in the case $S<1$) and the assumption $c\geq a$ (in the case $S>1$), we have 
$\mathbb{E}\left[\int_{0}^{{\tau}}c_{s}^{1-S}ds\right]<+\infty.$
 Therefore, for all $S > 0$, we obtain
 $$\mathbb{E}\left[((1-R)J({\tau},X_{\tau}^{0,x,\pi,c}))^{\frac{1}{\nu}}+\int_{0}^{{\tau}}c_{s}^{1-S}ds\right]<+\infty.$$
Finally, by Lemma \ref{lem:bsde-rs<1}, the BSDE (\ref{bsdegdpp}) admits a unique solution. \hfill$\Box$

\begin{lemma}\label{lemma3.6}
For   all $(t,x,\pi,c)\in\mathbb{R}_+\times[\frac{a}{r},+\infty)\times{\mathcal{A}}(x)$ 
and  $\varepsilon>0$ there exists an
admissible control $(\hat{\pi},\hat{c})\in{\mathcal{A}}(t,X_t^{0,x,\pi,c})$ such that
\begin{align}\label{3.15}
J(t,X_t^{0,x,\pi,c})\leq Y^{t,X_t^{0,x,\pi,c},\hat{\pi},\hat{c}}_t+\varepsilon, ~{\mathbb{P}}\mbox{-a.s.}.
\end{align}
\end{lemma}

\noindent\textbf{Proof.} Let $\{h^n\}$, $n\in N$,
be a dense subset of $\mathbb{R}_{+}$ with $h^1=a$, and define $B(h^n,\frac{1}{k})=[h^n,h^n+\frac{1}{k})$.
Set $B_{n,k}:=B(h^n,\frac{1}{k})\cap \mathbb{R}_{+}\setminus
\bigcup_{m<n}B(h^m,\frac{1}{k})$  and
$A_{n,k}:=\{\omega\in \Omega|X_t^{0,x,\pi,c}(\omega)\in
B_{n,k}\}$. Then $\cup_{n=1}^{\infty}A_{n,k}=\Omega$, $f_k:= \Sigma^{\infty}_{n=1}h^n1_{
	A_{n,k}}(\omega)\leq X_t^{0,x,\pi,c}$ and the
sequence $f^k(\omega)$ is ${\mathcal
	{F}}_{t}$-measurable and
converges to $X_t^{0,x,\pi,c}$ strongly and
uniformly. By the uniformly continuous of $J(x)$, and  Lemma \ref{25121201},
we can choose  $k$ sufficiently large such that 
$$
\left|J(t,f^k)-J(t,X_t^{0,x,\pi,c})\right|=\left|J(0,f^k)-J(0,X_t^{0,x,\pi,c})\right|\leq
\frac{\varepsilon}{2},\ {\mathbb{P}}\mbox{-a.s.}
$$

For every $h^n\in \mathbb{R}_+$,
we  choose an admissible control
$(\pi_n,c_n)\in \mathcal{A}(t,h^n)$ such that
$$
J(t,h^n)\leq
Y^{t,h^n, \pi_n,c_n}_t+\frac{\varepsilon}{2}, \quad  {\mathbb{P}}\mbox{-a.s.}.
$$
Then
$(\hat{\pi},\hat{c}):=\sum^{\infty}_{n=1}(\pi_n,c_n)1_{\{X_t^{0,x,\pi,c}\in
	B_{n,k}\}}\in \mathcal{A}(t,f_k)\subset \mathcal{A}(t,X_t^{0,x,\pi,c})$ and
\begin{align*}
Y^{t,X_t^{0,x,\pi,c},\hat{\pi},\hat{c}}_t&=
Y^{t,f^k,\hat{\pi},\hat{c}}_t=
\sum^{\infty}_{n=1}Y^{t,h^n,\pi_n,c_n}_t1_{\{X_t^{0,x,\pi,c}\in
		B_{n,k}\}}\\
&\geq-\frac{\varepsilon}{2}+\sum^{\infty}_{n=1}J(t,h^n)1_{\{X_t^{0,x,\pi,c}\in
		B_{n,k}\}}=-\frac{\varepsilon}{2}+J(t,f^k)
	\geq-\varepsilon+ J(t,X_t^{0,x,\pi,c}),  \ \  {\mathbb{P}}\mbox{-a.s.}.
\end{align*}
Thus, we have \eqref{3.15}. \hfill$\Box$

\subsection*{Proof of Theorem \ref{theoremddp111}}
 We first prove the equation (\ref{ddpG25}) for deterministic times $T > 0$. By the definition of $J(x)$ and Lemma  \ref{th:finte-to-infinte}, we have, for every deterministic time $T>0$,
$$
J(x)=J(0,x)=\mathop{\esssup}\limits_{(\pi,c)\in{\mathcal{A}}(x)}V^{0,x,\pi,c}_0=\mathop{\esssup}\limits_{(\pi,c)\in{\mathcal{A}}(x)}Y^{0,x,\pi,c}_0,
$$where $Y^{0,x,\pi,c}$ is the solution of BSDE \eqref{bsdetau}. 
By the uniqueness of solution to the BSDE \eqref{bsdegdpp}, 
it follows that $Y^{0,x,\pi,c}_0=G^{0,x,\pi,c}_{0,T}\left[Y^{0,x,\pi,c}_{T}\right]$. Moreover, from the uniqueness property of the BSDE \eqref{bsdetau}, we have $Y^{0,x,\pi,c}_{T}=Y^{T,X^{0,x,\pi,c}_{T},\pi,c}_{T}$. Then, we obtain
$$
J(x)
=\mathop{\esssup}\limits_{(\pi,c)\in{\mathcal{A}}(x)}G^{0,x,\pi,c}_{0,T}\left[Y^{T,X^{0,x,\pi,c}_{T},\pi,c}_{T}\right].
$$
For every $(\pi,c)\in \mathcal{A}(x)$, by Lemma 2.26 in \cite{fab1}, for $\mathbb{P}$-a.e. $\omega_0\in \Omega$,
the shifted control $(\pi^{\omega_0},c^{\omega_0})(\cdot):=(\pi,c)|_{[T,\infty)}$ is $\mathcal{F}_{\omega^T_0,s}$-progressively measurable and $ B^T_{\cdot}:=B_{\cdot}-B_T$ is a Wiener process on probability space $(\Omega,\mathcal{F}_{\omega_0}, \mathcal{F}_{\omega^T_0,s}, \mathbb{P}_{\omega_0})$, where  $\mathbb{P}_{\omega_0}=\mathbb{P}(\cdot|\mathcal{F}_T)(\omega_0)$ is the regular conditional probability, $\mathcal{F}_{\omega_0}$ is the augmentation of $\mathcal{F}$ by the  $P_{\omega_0}$ null sets, and $\mathcal{F}_{\omega^T_0,s}$ is the augmented filtration generated by $B^T$.

Since $(\pi,c)\in \mathcal{A}(x)$ implies  $X^{0,x,\pi,c}_s\geq 0$ for all $s\geq0$, we have $(\pi^{\omega_0},c^{\omega_0})\in \mathcal{A}(T,X^{0,x,\pi,c}_T(\omega_0))$ for $\mathbb{P}$-a.e. $\omega_0\in \Omega$. Therefore, using that for $\mathbb{P}$-a.e. $\omega_0$, $\mathbb{P}_{\omega_0}(\{\omega:X^{0,x,\pi,c}_T(\omega)=X^{0,x,\pi,c}_T(\omega_0)\})=1$, we obtain
\begin{align*}
&Y^{T,X^{0,x,\pi,c}_{T},\pi,c}_{T}(\omega_0)=Y^{T,X^{0,x,\pi,c}_{T},\pi^{\omega_0},c^{\omega_0}}_{T}(\omega_0)
=\mathbb{E}\left[Y^{T,X^{0,x,\pi,c}_{T},\pi^{\omega_0},c^{\omega_0}}_{T}|\mathcal{F}_T\right](\omega_0)\\
=&\mathbb{E}_{\omega_0}\left[Y^{T,X^{0,x,\pi,c}_{T},\pi^{\omega_0},c^{\omega_0}}_{T}\right]\leq \mathbb{E}_{\omega_0}\left[J(T,X^{0,x,\pi,c}_{T})\right]=J(T,X^{0,x,\pi,c}_{T}(\omega_0)).
\end{align*}
Then, by 
Lemma \ref{lem:bsde-rs<1}, we have
$$
J(x)=J(0,x)\leq \mathop{\esssup}\limits_{(\pi,c)\in{\mathcal{A}}(x)}G^{0,x,\pi,c}_{0,T}[J(T,X^{0,x,\pi,c}_{T})].
$$
Moreover, by Lemma \ref{25121201}, we have $J(T,X^{0,x,\pi,c}_{T})=J(T,y)|_{y=X^{0,x,\pi,c}_{T}}=J(y)|_{y=X^{0,x,\pi,c}_{T}}=J(X^{0,x,\pi,c}_{T}).$
Therefore, we obtain
\begin{align}
\label{eq:dpp-upper}
J(x)\leq\mathop{\esssup}\limits_{(\pi,c)\in{\mathcal{A}}(x)}G^{0,x,\pi,c}_{0,T}[J(X^{0,x,\pi,c}_{T})].
\end{align}

On the other hand, by Lemma \ref{lemma3.6}, for any
$\varepsilon>0$, there
is an admissible control $(\bar{\pi},\bar{c})\in{\mathcal{A}}(T,X^{0,x,\pi,c}_{T})$ such that
$$
J(T,X^{0,x,\pi,c}_{T})\leq Y^{T,X^{0,x,\pi,c}_{T},\bar{\pi},\bar{c}}_{T}+\varepsilon, \ \  \mathbb{P}\mbox{-a.s.}
$$
Then, by Lemma \ref{lem:bsde-rs<1},
 we have
\begin{align*}
\left|G^{0,x,\pi,c}_{0,T}\left[Y^{T,X^{0,x,\pi,c}_{T},\bar{\pi},\bar{c}}_{T}\right]
	-G^{0,x,\pi,c}_{0,T}[J(T, X^{0,x,\pi,c}_{T})]\right|\leq \varepsilon.
\end{align*}

    Define
$(\hat{\pi},\hat{c})=(\pi,c)1_{[0,T]}(\cdot)+(\bar{\pi},\bar{c})1_{(T,\infty)}(\cdot)
\in {\mathcal{A}}(x).
$ Therefore, 
\begin{align*}
J(x)=J(0,x)&\geq
G^{0,x,\hat{\pi},\hat{c}}_{0,T}\left[Y^{T,X^{0,x,\hat{\pi},\hat{c}}_{T},\hat{\pi},\hat{c}}_{T}\right]
=G^{0,x,\pi,c}_{0,T}\left[Y^{T,X^{0,x,\pi,c}_{T},\bar{\pi},\bar{c}}_{T}\right]\\
&\geq G^{0,x,\pi,c}_{0,T}[J(T,X^{0,x,\pi,c}_{T})]-\varepsilon
= G^{0,x,\pi,c}_{0,T}[J(X^{0,x,\pi,c}_{T})]-\varepsilon.
\end{align*}
By the arbitrariness of $\varepsilon$ in the last inequality, together with the upper bound in (\ref{eq:dpp-upper}), we obtain 
\begin{align*}
J(x)=\mathop{\esssup}\limits_{(\pi,c)\in{\mathcal{A}}(x)}G^{0,x,\pi,c}_{0,T}[J(X^{0,x,\pi,c}_{T})].
\end{align*}
Finally, by following the proof of Theorem 3.70 in \cite{fab1} (see pp. 241–244), we show that \eqref{ddpG25} holds for any bounded stopping time $\tau$. \hfill$\Box$

\section*{Appendix B: Proofs in Section \ref{sec:C2}}
\subsection*{Proof of Theorem \ref{prop:vis}}
First, we show $J$ is a viscosity subsolution of \eqref{eq:Case2-vis} by a contradiction argument. Let $\varphi\in {\cal{A}}^+(\hat{x},J)$, where $\hat{x}>\frac{a}{r}$. Without loss of generality, that the supremum in the definition of  ${\cal{A}}^+(\hat{x},J)$ is attained as a strict maximum.  
Assume on the contrary that 
$$-\delta\nu J(\hat{x})+{\mathbf{H}}(\hat{x},J(\hat{x}),\varphi'(\hat{x}),\varphi''(\hat{x}))<0.$$
By  the continuity of $J$ and ${\mathbf{H}}$, and the monotonic decreasing property of ${\mathbf{H}}$ in the second variable, there exist $h>0$ and $\varepsilon> 0 $ such that
\begin{align}\label{eq:xiajie-2e}
-\delta\nu z+{\mathbf{H}}(y,z,\varphi'(y),\varphi''(y))<-\varepsilon, \quad z\in [J(\hat{x})-h,+\infty),
\end{align}
for all $y\in B(\hat{x},h)=\{y:|\hat{x}-y|\leq h\}$, and 
\begin{align}\label{2606032}
(J-\varphi)(\hat{x}+h)\vee  (J-\varphi)(\hat{x}-h)\leq -2\varepsilon.
\end{align}
For fixed  $h>0$, define
$\tau^{c,h}:=\tau^{\pi,c,h}:=\inf\left\{t:|X^{\hat{x},\pi,c}_t-\hat{x}|\geq h
\right\},$ and let $ \quad (\pi,c)\in \mathcal{A}(\hat{x}).$
Then, for every $k>0$, by Theorem \ref{theoremddp111}, 
there exists $(\pi^{\varepsilon},c^{\varepsilon})\equiv (\pi^{\varepsilon,h,k},c^{\varepsilon,h,k})\in {\mathcal{A}}(\hat{x})$ such
that 
\begin{align}\label{2606031}
J(\hat{x})\leq\mathbb{E}\left[J(X^{\hat{x},\pi^{\varepsilon},c^{\varepsilon}}_{\tau^{c^\varepsilon,h}\wedge k})+\int_0^{\tau^{c^\varepsilon,h}\wedge k} [f(c^{\varepsilon}_s,Y^{c^{\varepsilon},J}_s)-\delta \nu Y^{c^{\varepsilon},J}_s]ds\right]+\varepsilon.
\end{align}
Applying  It\^{o}'s formula  to
 ${\varphi}(X^{\hat{x},\pi^{\varepsilon},{c}^{\varepsilon}}_s)$, we obtain
\begin{align*}                            \varphi(X^{\hat{x},\pi^{\varepsilon},c^{\varepsilon}}_{\tau^{c^\varepsilon,h}\wedge k})
=\varphi(\hat{x})+\int^{\tau^{c^\varepsilon,h}\wedge k}_{0} ({\cal{L}}\varphi)(X^{\hat{x},\pi^{\varepsilon},{c}^{\varepsilon}}_s,
{\pi}^{\varepsilon}_s, c^{\varepsilon}_s)ds
+\int^{\tau^{c^\varepsilon,h}\wedge k}_{0}
\pi_s^{\varepsilon}\sigma\partial_x{\varphi}(X^{\hat{x},\pi^{\varepsilon},{c}^{\varepsilon}}_s)dB_s,
\end{align*}
where, for any $( x,\pi,c)\in [\frac{a}{r},+\infty)\times \mathbb{R}\times \mathbb{R}_+$, 
\begin{align*}
({\cal{L}}{\varphi})(x,\pi,c) :=
\langle{\varphi}'(x),rx+(\mu-r)\pi-c\rangle+\frac{1}{2}\varphi''(x)\pi^2\sigma^2.
\end{align*}
Notice that, by the definition of $\mathbf{H}$ and  \eqref{eq:xiajie-2e},
\begin{align*}
\delta\nu z-\varepsilon&\geq \mathbf{H}(x,z,\varphi'(x),\varphi''(x))\geq({\cal{L}}{\varphi})(x,\pi,c)+\sup_{c\geq a}[f(c,z)-c\varphi'(x)]+c\varphi'(x)\nonumber\\
&\geq ({\cal{L}}{\varphi})(x,\pi,c)+f(c,z), \quad z\in [J(\hat{x}),+\infty),\ x\in B(\hat{x},h).
\end{align*}
It follows that  
\begin{align}\label{bsde4.210603}
J(\hat{x})=\varphi(\hat{x})&=\mathbb{E}\varphi(X^{\hat{x},\pi^{\varepsilon},c^{\varepsilon}}_{\tau^{c^\varepsilon,h}\wedge k})-
\mathbb{E}\int^{\tau^{c^\varepsilon,h}\wedge k}_{0} ({\cal{L}}\varphi)(X^{\hat{x},\pi^{\varepsilon},{c}^{\varepsilon}}_s, {\pi}^{\varepsilon}_s, c^{\varepsilon}_s)ds\nonumber\\
  &\geq    \mathbb{E}\varphi(X^{\hat{x},\pi^{\varepsilon},c^{\varepsilon}}_{\tau^{c^\varepsilon,h}\wedge k}) +\mathbb{E}  \int_0^{\tau^{c^\varepsilon,h}\wedge k} [f(c^{\varepsilon}_s,Y^{c^{\varepsilon},J}_s)-\delta \nu Y^{c^{\varepsilon},J}_s]ds+\varepsilon\mathbb{E}[\tau^{c^\varepsilon,h}\wedge k].
\end{align}
Combining \eqref{2606031} and \eqref{bsde4.210603} together, we obtain 
 \begin{align}\label{2606033}
\varepsilon k\mathbb{P}(\tau^{c^\varepsilon,h}>k)
&\leq\varepsilon k\mathbb{P}(\tau^{c^\varepsilon,h}>k)+\varepsilon \mathbb{E}[\tau^{c^\varepsilon,h}1_{\tau^{c^\varepsilon,h}\leq k}]\nonumber\\
&=\varepsilon\mathbb{E}[\tau^{c^\varepsilon,h}\wedge k]\leq \mathbb{E}[(J-\varphi)(X^{\hat{x},\pi^{\varepsilon},c^{\varepsilon}}_{\tau^{c^\varepsilon,h}\wedge k})]+\varepsilon\leq-2\varepsilon\mathbb{P}(\tau^{c^\varepsilon,h}\leq k)+\varepsilon, 
 \end{align}where the last inequality follows from  
 \eqref{2606032}. Dividing both sides of \eqref{2606033} by $\varepsilon k$, it implies 
\begin{align*}
\mathbb{P}(\tau^{c^\varepsilon,h}>k)
\leq \frac{1}{k}, \text{ and thus, } \quad  \mathbb{P}(\tau^{c^\varepsilon,h}\leq k)
\geq \frac{k-1}{k}.
 \end{align*}
Let $k\geq 4$. Then, by the last inequality and \eqref{2606033}, we obtain
$
 0\leq -2\varepsilon\mathbb{P}(\tau^{c^\varepsilon,h}\leq k)+\varepsilon\leq -\frac{\varepsilon}{2}.
$
Obviously, it is a contradiction.

Second, we show that $J$ is a viscosity supersolution of \eqref{eq:Case2-vis} and therefy finish the proof. 

Let  $\varphi\in {\cal{A}}^-(\hat{x},J)$ with
 $\hat{x}>\frac{a}{r}$.  For any given $(\bar{\pi},\bar{c})\in \mathbb{R}\times [a,+\infty)$, let $(\pi,c)\in\mathcal{A}(\hat{x})$
such that $(\pi_0,c_0)=(\bar{\pi},\bar{c})$ and are continuous on $t=0$. For example, we can let $(\pi_t,c_t)=(\bar{\pi},\bar{c})1_{[0,\tau)}+(0,a)1_{[\tau,\infty)}$, where $\tau=\inf\{t: X_t<\frac{1}{2}(\frac{a}{r}+\hat{x})\}$. Therefore, by Theorem \ref{theoremddp111}, we  obtain that  for any $h>0$,
 \begin{align*}
0=J(\hat{x})+\varphi (\hat{x})
&\geq\mathbb{E}\left[J(X^{\hat{x},\pi,c}_{\tau\wedge h})+\int_0^{\tau\wedge h}[f(c_s,Y^{c,J}_s)-\delta \nu Y^{c,J}_s]ds\right]+\varphi (\hat{x})\\
&\geq \mathbb{E}\left[-\varphi(X^{\hat{x},\pi,c}_{\tau\wedge h})
+\int_0^{\tau\wedge h}[f(c_s,Y^{c,J}_s)-\delta \nu Y^{c,J}_s]ds\right]+\varphi (\hat{x}).
\end{align*}
Applying  It\^{o} formula  to
${\varphi}(X^{\hat{x},\pi,{c}}_s)$, we have
\begin{align*}
 0 \geq\frac{1}{h}\mathbb{E}\left[\int_0^{\tau\wedge h}[-({\cal{L}}\varphi)(X^{\hat{x},\pi,{c}}_s,
{\pi}_s, c_s)+f(c_s,Y^{c,J}_s)-\delta \nu Y^{c,J}_s]ds\right].
\end{align*}
Let $h\rightarrow0$  we obtain 
$$-({\cal{L}}\varphi)(\hat{x},\bar{\pi}, \bar{c})+f(\bar{c},J(\hat{x}))-\delta \nu J(\hat{x})\leq0.$$
By the arbitrariness of $(\bar{\pi},\bar{c})$,
$$-\delta\nu J(\hat{x})+{\mathbf{H}}(\hat{x},J(\hat{x}),-\varphi'(\hat{x}),-\varphi''(\hat{x}))\leq0.$$It completes the proof. 
\hfill$\Box$

\subsection*{Proof of Proposition \ref{prop:recursive-C1}}
Since the function $J$ is increasing and concave by Propositions \ref{pro:ez-increase}  and \ref{pro:ez-concave}, we define its right and left derivatives for $x > \frac{a}{r}$ by
\begin{equation*}
J'_{\pm}(x)=\lim_{h\rightarrow 0+} {J(x\pm h)-J(x)\over \pm h} \geq 0.
\end{equation*} 
To establish that $J$ is $C^1$,  it suffices to show that $J$ is differentiable, i.e., $J'_{+}(x)= J'_{-}(x)$ for all $x>a/r$.

We prove the result by contradiction. Assume that $J'_{+}(x_0)< J'_{-}(x_0)$ for some $x_0 \in ({a\over r},+\infty)$. Let $\beta$ be a constant  such that $J'_{+}(x_0)<\beta< J'_{-}(x_0)$. For each $m>0$, we define a smooth function, 
\begin{equation}
\label{eq: def of phi-recursive}
\phi(x)=J(x_0)+\beta(x-x_0)-m(x-x_0)^2, ~~x>a/r.
\end{equation}
It is clear that $\phi(x_0)=J(x_0)$, $\phi'(x_0)=\beta$ and $\phi''(x_0)=-2m$. By the concavity of $J$, for $0<x_0-x<{1\over m}(J'_{-}(x_0)-\beta)$, we have
$
J(x) \leq  J(x_0)+J'_{-}(x_0)(x-x_0). 
$
Substituting the expression for $\phi$ from equation (\ref{eq: def of phi-recursive}) into the inequality above, we obtain
\begin{align}
\label{eq:test-left-X-recursive}
J(x)\leq \phi(x)+(J'_{-}(x_0)-\beta)(x-x_0)+m(x-x_0)^2<\phi(x).
\end{align}
Similarly, for $0<x-x_0<{1\over m} (\beta-J'_{+}(x_0))$, we have
\begin{align}
\label{eq:test-right-X-recursive}
J(x)&\leq  J(x_0)+J'_{+}(x_0)(x-x_0)\\ \nonumber
&=\phi(x)+(J'_{+}(x_0)-\beta)(x-x_0)+m(x-x_0)^2< \phi(x).
\end{align}
Equations (\ref{eq:test-left-X-recursive}) and (\ref{eq:test-right-X-recursive}) together imply that $J(x)<\phi(x)$ in a small neighborhood of $x_0$. Consequently, $\phi(x)$ is a test function at $x=x_0$.

By Theorem \ref{prop:vis}, $J$ is a viscosity solution of  
\begin{equation}
\label{eq:visocisity}
 -\kappa {(J')^2\over J''}+G(J,J')+rxJ'-\delta\nu J=0,
\end{equation}
where $ G(v,d)=\sup_{c\geq a}[f(c,v)-cd]$ for $v\in \mathbb{V}$ and $d>0$. Using the test function $\phi(x)$ at $x_0$, we obtain
\begin{align*}
0&\geq \kappa{(\phi'(x_0))^2\over \phi''(x_0)}- G(\phi(x_0),\phi'(x_0)) -rx_0  \phi'(x_0)+\delta\nu \phi(x_0) \\\nonumber
&=-{\kappa \beta^2\over 2 m}- G(J(x_0),\beta) -rx_0 \beta+\delta\nu J(x_0). 
\end{align*}
Sending $m\rightarrow+\infty$, the last inequality implies that
\begin{align}
\label{eq: g1-X-recursive}
0 \geq  - G(J(x_0),\beta) -rx_0 \beta+\delta \nu J(x_0),~~  \beta\in (J'_{+}(x_0),J'_{-}(x_0)).
\end{align}

On the other hand, since $J(\cdot)$ is concave, by Alexandrov theorem, it is twice differentiable almost everywhere. Then, there exists an increase sequence $\{x_n\}_{n\geq 0}$ such that $x_n \rightarrow x_0$, and $J(\cdot)$ is $C^2$ at all $x_n$. By (\ref{eq:visocisity}), we have
\begin{align*}
0&= \kappa {(J'(x_n))^2\over J''(x_n)}-G(J(x_n),J'(x_n))-rx_nJ'(x_n)+\delta\nu J (x_n)  \\
&\leq  -G(J(x_n),J'(x_n))-rx_nJ'(x_n)+\delta\nu J (x_n).
\end{align*}
Taking the limit as $x_n \uparrow x_0$, since $G$ is continuous in $v$ and $d$, we obtain
\begin{align}
\label{eq:g2-X-recursive}
0\leq -G(J(x_0),J'_{-}(x_0))-rx_0J'_{-}(x_0)+\delta\nu J (x_0).
\end{align}
Similarly, by choosing $x_n \downarrow x_0$, we get
\begin{align}
\label{eq:g3-X-recursive}
0\leq -G(J(x_0),J'_{+}(x_0))-rx_0J'_{+}(x_0)+\delta\nu J (x_0).
\end{align}
Define $\varphi_0(d)=G(J(x_0),d)+rx_0d-\delta\nu J(x_0)$. It is clear that $\varphi_0(\cdot)$ is a convex function.  By \eqref{eq:g2-X-recursive} and \eqref{eq:g3-X-recursive}, it implies that
\begin{align*}
\varphi_0(d) \leq \max (\varphi_0(J'_{-}(x_0)),\varphi_0(J'_{+}(x_0)) )\leq 0,   \quad \forall d\in (J'_{+}(x_0),J'_{-}(x_0)).
\end{align*}
However, \eqref{eq: g1-X-recursive} implies that
$\varphi_0(d)  \geq 0,  \quad \forall d\in (J'_{+}(x_0),J'_{-}(x_0)).$
It thus follows that
$\varphi_0(d)\equiv 0$ on $d \in (J'_{+}(x_0),J'_{-}(x_0))$.  
However, by a straightforward calculation,  
\begin{align*}
\varphi_0(d)=
\begin{cases}
{S\over 1-S}[((1-R)J(x_0))]^{\rho\over S}d^{1-{1\over S}}+r x_0d-\delta\nu J(x_0),  & [(1-R)J(x_0)]^{\rho/S} d^{-1/S}>a,  \\
f(a,J(x_0))-ad+r x_0d-\delta\nu J(x_0),  & [(1-R)J(x_0)]^{\rho/S} d^{-1/S}\leq a. 
\end{cases}
\end{align*}
Therefore, it is impossible that $\varphi_0(d)\equiv 0$ on $d \in (J'_{+}(x_0),J'_{-}(x_0))$ unless $J'_{+}(x_0) = J'_{-}(x_0)$. 
Thus, we conclude that the value function $J$ is $C^1$ for all $x > \frac{a}{r}$. \hfill$\Box$

\subsection*{Proof of Proposition \ref{prop:strictly increasing}}

(i) By Proposition \ref{prop:recursive-C1},  $J$ is $C^1$ in $({a\over r},+\infty)$. Suppose, for contradiction, that there exists $x_0>\frac{a}{r}$ such that $J'(x_0)=0$. Given that $J$ is both increasing and concave by Proposition \ref{pro:ez-increase} and Proposition \ref{pro:ez-concave}, it follows that $J'(x)=0,  \forall  x\in [x_0,+\infty)$. However, this contradicts the inequality in (\ref{eq:bound}).

(ii) We  prove that for the point such that $J''(x_0)$ exists, it must hold that $J''(x_0)\neq 0$. Assume not, for each $m>0$, we define a smooth function
\begin{equation*}
\phi(x)=J(x_0)+J'(x_0)(x-x_0)-m(x-x_0)^2, ~~x>a/r.
\end{equation*}
Since $J''(x_0)=0$, for any $m>0$, there exists a neighborhood of $x_0$ over which $\phi(x)<J(x)$. We now apply the viscosity supersolution property and get
\begin{align*}
0&\leq \kappa {(J'(x_0))^2\over J''(x_0)}-G(J(x_0),J'(x_0))-rx_0J'(x_0)+\delta\nu J (x_0)  \\
&= -\kappa {{(J'(x_0))^2\over 2m}} -G(J(x_0),J'(x_0))-rx_0J'(x_0)+\delta\nu J (x_0).
\end{align*}
Since the terms $G(J(x_0),J'(x_0)),rx_0J'(x_0),\delta\nu J (x_0)$ are all finite,  sending $m \rightarrow 0$ leads to the desired contradiction.\hfill$\Box$


\subsection*{Proof of Theorem \ref{thm:recursive}}

Since $J(\cdot)$ is increasing, we have $\lim_{x\rightarrow +\infty}J'(x)\geq 0$. Moreover, by the concavity of $J$ and \eqref{eq:bound}, we have
\begin{equation*}
J'(x)\leq {J(x)-J({a\over r})\over x-{a\over r}}\leq {\eta^{-\nu S}{x^{1-R}\over 1-R}-J(\frac{a}{r})\over x-{a\over r}}, ~~x>\frac{a}{r}.
\end{equation*}
Then $\lim_{x\rightarrow +\infty}J'(x)=0$ for any $R > 0$ (we use L' Hospital rule for $0<R<1$).

We define the dual transformation
\begin{equation}
\label{eq: def-dual-X}
v(y):= \max_{x>{a\over r}}(J(x)-xy),  \quad 0=\lim_{x\rightarrow +\infty}J'(x)<y< J'({a\over r}) := \lim_{x \downarrow \frac{a}{r}}J'(x).
\end{equation} 
  Then $v(\cdot)$ is a decreasing and convex function on $(0,J'({a\over r}))$. By Proposition \ref{prop:strictly increasing}, $J'(\cdot)$ is strictly decreasing. We denote the inverse function of $J'(x)=y$ by $
I(y)=x$. The function $I(\cdot)$ is strictly decreasing and maps $(0,J'({a\over r}))$ to $({a\over r},+\infty)$.


  We next show the value function $J$ is $C^2$ when $x\in \mathcal U$ and $x\in \mathcal B$. By (\ref{eq: def-dual-X}), we obtain
\begin{equation}
\label{eq:dual-plug-X}
v(y)=[J(x)-xJ'(x)]|_{x=I(y)}=J(I(y))-yI(y),~~y\in(0,J'({a\over r})).
\end{equation}

We consider the following ODE of the function $v(\cdot)$ on $(0,J'({a\over r}))$, 
\begin{equation}
\label{eq:transODE-X}
\delta\nu(v(y)-yv'(y))=\kappa y^2v''(y)+ G(v(y)-yv'(y),y) - r v'(y) y.
\end{equation}
Equation (\ref{eq:transODE-X}) is a quasilinear ODE that degenerates only at $y=0$. The coefficient of $v''(y)$ in the above equation is $\kappa y^2$, which is nonzero when $y>0$. By a similar argument in Theorem 7.1 in Xu and Yi (2016), we have
\begin{equation*}
v\in C^2(0, J'({a\over r})).
\end{equation*}
Differentiating (\ref{eq:dual-plug-X}) once and twice,  we get
\begin{equation}
\label{eq:dual-first-X}
v'(y)=J'(I(y))I'(y)-yI'(y)-I(y)=-I(y),
\end{equation}
and
\begin{equation}
\label{eq:dual-second-X}
v''(y)=-I'(y)=-{1\over J''(I(y))}.
\end{equation}
By Proposition \ref{prop:strictly increasing} (ii), $J''(I(y))\neq 0$. Combining (\ref{eq:dual-plug-X}) and (\ref{eq:dual-first-X}), we get
\begin{equation}
\label{eq:V-dual-X}
J(I(y))=v(y)-yv'(y).
\end{equation}
Now, substituting (\ref{eq:dual-first-X})-(\ref{eq:V-dual-X}) into (\ref{eq:visocisity}) and recalling that $G(v,d)=\sup\limits_{c\geq a}[f(c,v)-cd]$,  we obtain the same equation (\ref{eq:transODE-X}). By (\ref{eq:dual-second-X}),  we establish that  $J(\cdot)$ is $C^2$ when $x\in \mathcal U$ and $x\in \mathcal{B}$.

Finally, we show that $J(x)$ is $C^2$ at the boundary point in  $cl({\mathcal U}) \bigcap cl({\mathcal B})$. It has already been shown that the value function $J(x)$ is $C^2$ when $x\in \mathcal U$ and $x\in \mathcal B$. Moreover, 
\begin{equation}\label{eq:J in U}
-\kappa {(J')^2\over J''}+{S\over 1-S}((1-R)J)^{\rho\over S}(J')^{1-{1\over S}}+rxJ'-\delta\nu J=0, ~~ \ \ x\in\mathcal{U};
\end{equation}
and  
\begin{equation}\label{eq:J in B}
-\kappa {(J')^2\over J''}+f(a,J)-aJ'+rxJ'-\delta\nu J=0, ~\ \ x\in \mathcal{B}.
\end{equation}

Assume $x$ is  the  connection point in $cl({\mathcal U}) \bigcap cl({\mathcal B})$.   Without loss of generality,  we assume that the left neighborhood of $x$ is $\mathcal B$ and right neighborhood of $x$  is $\mathcal U$.  Then, we obtain 
\begin{equation}
\label{eq:ODE-left}
\delta\nu J(x-)=-\kappa {(J'(x-))^2\over J''(x-)}+f(a,J(x-))-aJ'(x-)+rxJ'(x-),
\end{equation}
and 
\begin{align}
\label{eq:ODE-right}
\delta\nu J(x+)=-\kappa {(J'(x+))^2\over J''(x+)}+{S\over 1-S}((1-R)J(x+))^{\rho\over S}J'(x+)^{1-{1\over S}}+rxJ'(x+).
\end{align}
Since $J$ is continuous on $({a\over r},+\infty)$, we have $
J(x+)=J(x-).$ Moreover, 
by Proposition \ref{prop:recursive-C1}, we obtain
\begin{equation}
\label{eq: der-match}
J'(x+)=J'(x-).
\end{equation}
Furthermore, since $\mathcal U$  is the right neighborhood of $x$, then for the optimal $c^*$, we have
\begin{equation}
\label{eq:opcontrol-unconstrained}
 c^*(x+)=((1-R)J(x+))^{\rho/S} J'(x+)^{-1/S}=a.
\end{equation}
Substituting (\ref{eq: der-match}) and (\ref{eq:opcontrol-unconstrained}) into (\ref{eq:ODE-left}) and (\ref{eq:ODE-right}),  we finally get that
$
J''(x+)=J''(x-).
$
Since $x$ is arbitrary connection point in $cl({\mathcal U}) \bigcap cl({\mathcal B})$,  we have demonstrated that $J$ is $C^2$ across the boundary between the two regions. The proof is complete. \hfill$\Box$

\renewcommand {\theequation}{B-\arabic{equation}} \setcounter
{equation}{0}
\renewcommand {\thelemma}{B.\arabic{lemma}} \setcounter
{theorem}{0}
\setcounter{equation}{0}
\section*{Appendix C:  Proofs in Section \ref{sec:verification}}

\subsection*{Proof of Proposition \ref{pro:J-daoshu}}
Let $V_t^{a}$ denote the unique utility process when $x = \frac{a}{r}$. By Proposition  \ref{pro:ez-increase}, we have $ V_t^a=\delta^{-\nu}\frac{a^{1-R}}{1-R}e^{-\delta\nu t}, \forall t\geq0.$ For the subsequent proof, we construct a sequence of admissible strategies $(\pi,c)$ as follows. For any $\xi>0$ and $\zeta\in\mathbb{R}$, we define 
$$c_t=a+\xi(X_t-\frac{a}{r})\text{~~~and~~~} \pi_tX_t=\zeta(X_t-\frac{a}{r}),$$where $X_t$ denotes the corresponding wealth process. 
Clearly, 
$
X_t  =  \frac{a}{r} + (x - \frac{a}{r}) e^{ l(\zeta, \xi) t + \sigma \zeta B_t },
$
where $l(\zeta, \xi) = r + \zeta(\mu-r) - \xi - \frac{1}{2} \zeta^2 \sigma^2$. For any $x > \frac{a}{r}$, it follows that $X_t > \frac{a}{r}$, and thus $(\pi, c)\in\mathcal{A}(x)$.  
We denote this choice of consumption rate as $c_t:=c_t(\zeta, \xi)$. The corresponding stochastic utility process is then written as $V^{c}:=V^{c(\zeta, \xi)}$. Moreover, 
we obtain
\begin{align}\label{eq:e-c-a}
\left(\mathbb{E}\big[ (c_s-a)^{\frac{1}{1-S}}\big]\right)^{1-S}
=\xi (x-\frac{a}{r})e^{l(\zeta, \xi) s+\frac{\sigma^2 \zeta^2}{2(1-S)}s}.\end{align}

 We first prove the inequality (\ref{eq:jlimit}) below for $R > 1$. Notice that $V^c\geq V^a$, and  $V^c$ is uniformly integrable. 
Moreover, $f(c,v)$ is joint concave in $(c,v)$ in this case. Then, we have
\begin{align*}
    V^c_t-V^a_t&=\mathbb{E}\left[\int_t^{\infty}e^{-\delta s}(f(c_s,V_s^c)-f(a,V_s^a))ds~\big|\mathcal{F}_t~\right]\\
&\geq \mathbb{E}\left[\int_t^{\infty}e^{-\delta s}\left(\frac{\partial f}{\partial c}(c_s,V_s^c)(c_s-a)+\frac{\partial f}{\partial v}(a,V_s^c)(V_s^c-V_s^a)\right)ds~\big|\mathcal{F}_t~\right]\\
&\geq \mathbb{E}\left[\int_t^{\infty}e^{-\delta s}\left(\frac{\partial f}{\partial c}(c_s,V_s^a)(c_s-a)+\frac{\partial f}{\partial v}(a,V_s^c)(V_s^c-V_s^a)\right)ds~\big|\mathcal{F}_t~\right],
\end{align*}
where the last inequality follows from the fact that $\frac{\partial^2 f(c,v)}{\partial c \partial v}>0$ and $V_s^c\geq V^a_s$.  
Applying the Gronwall's inequality and the reverse H\"{o}lder inequality, we deduce that
\begin{align*}
    V^c_t-V^a_t&\geq \mathbb{E}\left[\int_t^{\infty}e^{-\delta s}\frac{\partial f}{\partial c}(c_s,V_s^a)e^{\int_t^s e^{-\delta u}\frac{\partial f}{\partial v}(a,V_u^c)du}(c_s-a)
    ds~\big|\mathcal{F}_t~\right]\\
    &\geq \int_t^{\infty}e^{-\delta s}\left(\mathbb{E}\big[ \left( \frac{\partial f}{\partial c}(c_s,V_s^a)e^{\int_t^s e^{-\delta u}\frac{\partial f}{\partial v}(a,V_u^c)du}\right)^{\frac{1}{S}}~\big|\mathcal{F}_t~\big]\right)^S \left(\mathbb{E}\big[ (c_s-a)^{\frac{1}{1-S}}~\big|\mathcal{F}_t~\big]\right)^{1-S}ds.
\end{align*}
Taking $t=0$, we obtain
\begin{align*}
    V^c_0-V^a_0 
\geq \int_0^{\infty}e^{-\delta s}\left(\mathbb{E}\big[ \left( \frac{\partial f}{\partial c}(c_s,V_s^a)e^{\int_0^s e^{-\delta u}\frac{\partial f}{\partial v}(a,V_u^c)du}\right)^{\frac{1}{S}}\big]\right)^S \left(\mathbb{E}\big[ (c_s-a)^{\frac{1}{1-S}}\big]\right)^{1-S}ds.
\end{align*}
Therefore, by equation (\ref{eq:e-c-a}), we conclude that
\begin{align*}
    \frac{J(x)-J(\frac{a}{r})}{x-\frac{a}{r}}&\geq \frac{V_0^c-V_0^a}{x-\frac{a}{r}} \geq \xi\int_0^{\infty}e^{-\delta s}\left(\mathbb{E}\big[ \left( \frac{\partial f}{\partial c}(c_s,V_s^a)e^{\int_0^s e^{-\delta u}\frac{\partial f}{\partial v}(a,V_u^c)du}\right)^{\frac{1}{S}}\big]\right)^S e^{l(\zeta, \xi) s+\frac{\sigma^2 \zeta^2}{2(1-S)}s}ds. 
\end{align*}
Again, since $f(\cdot,v)$ is concave and $c_t\geq a$ and $\frac{\partial f}{\partial v}<0$, we have
$$\left(\mathbb{E}\big[ \left( \frac{\partial f}{\partial c}(c_s,V_s^a)e^{\int_0^s e^{-\delta u}\frac{\partial f}{\partial v}(a,V_u^c)du}\right)^{\frac{1}{S}}\big]\right)^S\leq \frac{\partial f}{\partial c}(a,V_s^a)=a^{-S}((1-R)V_s^a)^{\rho}.$$ 
Moreover, when $x\downarrow \frac{a}{r}$, it follows that $X_s\downarrow \frac{a}{r}$, which implies that $c_s\downarrow a$ and  $V_s^c\downarrow V_s^a$ (Theorem 6.5 in \cite{HHJ23b}). By the dominance convergence theorem,  for any $\xi>0$ and $\zeta\in\mathbb{R}$, this implies that,
\begin{align}\label{eq:jlimit}
    \lim_{x\downarrow \frac{a}{r}}\frac{J(x)-J(\frac{a}{r})}{x-\frac{a}{r}}\geq \xi\int_0^{\infty}e^{-\delta s} \frac{\partial f}{\partial c}(a,V_s^a)e^{\int_0^s e^{-\delta u}\frac{\partial f}{\partial v}(a,V_u^a)du} e^{l(\zeta, \xi) s+\frac{\sigma^2 \zeta^2}{2(1-S)}s}ds. 
\end{align}

 We next prove the inequality (\ref{eq:jlimit}) above for $0 < R< 1$.
In this case, 
 the function $f(c,v)$ is concave in $c$ and convex in $v$. Then, we obtain
\begin{align*}
    V^{c}_t-V^a_t&=\mathbb{E}\left[\int_t^{\infty}e^{-\delta s}(f(c_s,V_s^{c})-f(c_s,V_s^a)+f(c_s,V_s^a)-f(a,V_s^a))ds~\big|\mathcal{F}_t~\right]\\
&\geq \mathbb{E}\left[\int_t^{\infty}e^{-\delta s}\left(\frac{\partial f}{\partial v}(c_s,V_s^a)(V_s^{c}-V_s^a)+\frac{\partial f}{\partial c}(c_s,V_s^a)(c_s-a)\right)ds~\big|\mathcal{F}_t~\right].
\end{align*}
Given the facts that $\frac{\partial f}{\partial v}<0$, $\frac{\partial^2 f}{\partial v\partial c}<0$ and $0<\frac{\partial f}{\partial c}(c_s,V_s^a)$, by applying the Gronwall's inequality and reverse H\"{o}lder inequality, we obtain
\begin{align*}
    V^{c}_0-V^a_0&\geq \mathbb{E}\left[\int_0^{\infty}e^{-\delta s}\frac{\partial f}{\partial c}(c_s,V_s^a)e^{\int_0^s e^{-\delta u}\frac{\partial f}{\partial v}(c_u,V_u^a)du}(c_s-a)
    ds\right]\\
    &\geq \int_0^{\infty}e^{-\delta s}\left(\mathbb{E}\big[ \left( \frac{\partial f}{\partial c}(c_s,V_s^a)e^{\int_t^s e^{-\delta u}\frac{\partial f}{\partial v}(c_u,V_u^a)du}\right)^{\frac{1}{S}}\big]\right)^S \left(\mathbb{E}\big[ (c_s-a)^{\frac{1}{1-S}}\big]\right)^{1-S}ds.
\end{align*}
Since $f(\cdot,v)$ is concave and $c_t>a$ and $\frac{\partial f}{\partial v}<0$ and $\frac{\partial^2 f}{\partial v\partial c}<0$,  we have
$$\left(\mathbb{E}\big[ \left( \frac{\partial f}{\partial c}(c_s,V_s^a)e^{\int_0^s e^{-\delta u}\frac{\partial f}{\partial v}(c_u,V_u^a)du}\right)^{\frac{1}{S}}\big]\right)^S\leq \frac{\partial f}{\partial c}(a,V_s^a)=a^{-S}((1-R)V_s^a)^{\rho}.$$
By equation \eqref{eq:e-c-a} and letting $x\downarrow \frac{a}{r}$, we have $c_t\downarrow a$, and consequently, \eqref{eq:jlimit} also holds for $0<R<1$ by dominance convergence theorem.  

We now prove that $J^{\prime}(\frac{a}{r})=+\infty$ by using equation \eqref{eq:jlimit}. Since $(1-R)V_s^a=\delta^{-\nu}a^{1-R}e^{-\delta\nu s}$, we obtain
\begin{align*}
  \frac{\partial f}{\partial c}(a,V_s^a)e^{\int_0^s e^{-\delta u}\frac{\partial f}{\partial v}(a,V_u^a)du} &=a^{-S}((1-R)V_s^a)^{\rho}e^{\int_0^s e^{-\delta u}\rho \nu a^{1-S}((1-R)V^a_u)^{\rho-1} du} = a^{-R}\delta^{-\nu\rho}.
\end{align*}
Then, \eqref{eq:jlimit} implies that
\begin{align}
\lim_{x\rightarrow \frac{a}{r}}\frac{J(x) - J(a/r)}{x - a/r} &\ge a^{-R}  \delta^{-\nu\rho} \xi  \int_{0}^{\infty} e^{-\delta t} e^{l(\zeta, \xi)t+\frac{\sigma^{2}\zeta^{2}}{2(1-S)} t }dt\nonumber\\
&=a^{-R} \delta^{-\nu\rho}  \frac{\xi}{\delta-l(\zeta,\xi)-\frac{\sigma^{2}\zeta^{2}}{2(1-S)} }. \label{eq:xizeta}
\end{align}
To finish the proof, it suffices to select $(\zeta,\xi)$ such that the term $\delta-l(\zeta,\xi)-\frac{\sigma^{2}\zeta^{2}}{2(1-S)} = \delta-r+\xi-\zeta(\mu-r)-\frac{S}{2(1-S)}\sigma^{2}\zeta^{2} \downarrow 0$, the right side of \eqref{eq:xizeta} diverges to infinity.

To this end, we divide the analysis into two separate cases. For $R>1$ (hence $S>1$), under the condition $r-\delta+{S-1\over S}\kappa>0$, we may choose $\zeta={(S-1)(\mu-r)\over S\sigma^2}$ and $\xi \downarrow r-\delta+{(S-1)\over S}\kappa>0$. This choice leads to $\delta-l(\zeta,\xi)-\frac{\sigma^{2}\zeta^{2}}{2(1-S)} \downarrow 0$. For $0<R<1$, (hence $0<S<1$), we choose $\xi$ large enough such that the discriminant $2\sigma^{2}\big[\kappa+\frac{S}{1-S}(\delta-r+\xi)\big]\geq 0$.
Let $\zeta\downarrow-{(\mu-r)+\sqrt{2\sigma^{2}[\kappa+S(\delta-r+\xi)]/(1-S)}\over {\sigma^2S/(1-S)}}$. This choice leads to  
$\delta-l(\zeta,\xi)-\frac{\sigma^{2}\zeta^{2}}{2(1-S)} \downarrow 0$. This concludes the proof. \hfill$\Box$


\subsection*{Proof of Proposition \ref{pro-jprime-guji}}
By \eqref{eq:hjb},  and differentiating both sides of equations \eqref{eq:J in U} -  \eqref{eq:J in B}, we obtain that 
\begin{align*}
  \left(-\kappa {(J')^2\over J''}+rxJ'-\delta\nu J\right)'&=((1-R)J)^{\rho/S}(J')^{-1/S}\left[J''-\frac{\rho\nu (J')^2}{(1-R)J}\right], ~x\in\mathcal{U}; \\
  \left(-\kappa {(J')^2\over J''}+rxJ'-\delta\nu J\right)'&=-\rho\nu a^{1-S}((1-R)J)^{\rho-1}J'+aJ'', ~x\in\mathcal{B}. 
\end{align*}Recalling the definitions of $\mathcal{U}$ and $\mathcal{B}$ in \eqref{eq:U} and \eqref{eq:B}, one can easily verify that when $((1-R)J)^{\rho/S}(J')^{-1/S}\rightarrow a$, then the right hand sides of the above two equations are equal. Hence, it implies $J\in C^3$.  

Since $J'>0$, $J''<0$, and assumption (A4) holds, for any $y>0$, we have 
\begin{align*}
 J'(x)\geq \frac{J(x+y)-J(x)}{y}&\geq \frac{\frac{\delta^{-\nu}r^{1-R} }{1-R}(x+y)^{1-R}-\frac{\eta^{-\nu S}}{1-R}x^{1-R}}{y}\\
 &\geq \frac{1}{\Delta}\left(\frac{\delta^{-\nu}r^{1-R} }{1-R}(1+\Delta)^{1-R}-\frac{\eta^{-\nu S}}{1-R}\right)x^{-R},
\end{align*}where the last inequality follows by taking $y=\Delta x$. Now choose $\Delta > 0$ sufficiently large so that  
$\frac{\delta^{-\nu}r^{1-R}(1+\Delta)^{1-R} }{1-R}>\frac{\eta^{-\nu S}}{1-R}$. It follows that there exists a constant $k_1>0$ such that  $J'(x)\geq k_1 x^{-R}$ for all $x>a/r$.

On the other hand, for any $x>a/r$ and any $y$ satisfying $0<y<x-\frac{a}{r}$, we have
\begin{align*}
 J'(x)\leq \frac{J(x)-J(x-y)}{y}&\leq \frac{\frac{\eta^{-\nu S}}{1-R}x^{1-R}-\frac{\delta^{-\nu}r^{1-R} }{1-R}(x-y)^{1-R}}{y}\\
 &\leq \frac{1}{1-\Delta}\left(\frac{\eta^{-\nu S}}{1-R}-\Delta^{1-R}\frac{\delta^{-\nu}r^{1-R} }{1-R}\right)x^{-R},
\end{align*}where the last inequality follows by taking $y=(1-\Delta)x \in (0, x - \frac{a}{r})$. Now choose  $\Delta \in (0,1)$  sufficiently small so that $\frac{\eta^{-\nu S}}{1-R} > \Delta^{1-R}\frac{\delta^{-\nu}r^{1-R} }{1-R} $. Let 
$x_0:=\frac{a}{\Delta r}>\frac{a}{r}$. Then, for any $x > x_0$, we have $y = (1-\Delta)x < x - \frac{a}{r}$, and the above inequality implies that  there exists a constant $k_2 > 0$ such that $J'(x)\leq k_2 x^{-R}$ for all $x>x_0>a/r$.
\hfill$\Box$

\subsection*{Proof of Proposition \ref{pro: SDE}}
By (A1) and (A3), $J'(\frac{a}{r})=+\infty$ and $J'>0$. Therefore,  showing that $X_t>\frac{a}{r}$ is equivalent to prove that 
$J'(X_t)$ is almost surely finite for all $t\geq0$. By Proposition \ref{pro-jprime-guji}, we have $J\in C^3$. Using It\^{o}'s formula to 
$J'(X_t)$, we obtain
\begin{align*}
    dJ'(X_t)
    &=\left\{J''(X_t)rX_t-2\kappa J'(X_t)-J''(X_t)c(X_t)+\kappa \frac{J'''(X_t)(J'(X_t))^2}{(J''(X_t))^2} \right\}dt -J'(X_t)\frac{\mu-r}{\sigma}dB_t.
\end{align*}
  By (A2) or \eqref{eq:visocisity}, we obtain 
$$-\kappa {(J')^2\over J''}+rxJ'-\delta\nu J+G(J,J')=0,$$
where $ G(J,J')=\sup_{c\geq a}[f(c,J)-cJ']$. Differentiating both sides with respect to $x$, we get
$$\kappa \frac{J'''(J')^2}{(J'')^2}-2\kappa J'+rxJ''+(r-\delta\nu)J'+\frac{\partial G}{\partial J}J'+ \frac{\partial G}{\partial J'}J''=0.$$By the envelope theorem, $\frac{\partial G}{\partial J}=\frac{\partial f}{\partial v}(c^*, J)=f_v(c^*,J)$ and $\frac{\partial G}{\partial J'}=-c^*$. Therefore, we have
$$\kappa \frac{J'''(J')^2}{(J'')^2}-2\kappa J'+rxJ''+(r-\delta\nu)J'+ J'\frac{\partial f}{\partial v}(c^*, J)-c^*J''=0,$$where
$c^*$ is the feedback strategy.
Substituting the last equation into the drift term of $J'(X_t)$ yields  
\begin{align}\label{eq:fankui-Z}
 dM_t=(\delta\nu-r-f_v(c(X_t),J(X_t)))M_tdt
    -\frac{\mu-r}{\sigma}M_tdB_t, ~~Z_0=J'(x),
\end{align}where $M_t:=J'(X_t)$, $X_t=(J')^{-1}(M_t)$ and 
$f_v(c(X_t),J(X_t))=\rho\nu (c(X_t))^{1-S}((1-R)J(X_t))^{\rho-1}$.

Next, we show that SDE \eqref{eq:fankui-Z} has a unique positive solution $M$.  In fact, it can be easily verified that  
$c(x)=\max\Big\{a, (J'(x))^{-1/S} \big((1-R)J\big)^{\frac{\rho}{S}}I_{x>\frac{a}{r}}\Big\}$ is continuous at $x=\frac{a}{r}$ by (A3). Hence, $f_v(c(x),J(x))$ is continuous at $[\frac{a}{r},+\infty)$. To show the boundedness of $f_{v}$, it  suffices to verify the bounded behaviors when $x\downarrow \frac{a}{r}$ and $x\uparrow+\infty$. 
Indeed, 
\begin{align*}
    \lim_{x\downarrow a/r}f_v(c(x),J(x))&=\rho\nu a^{1-S}((1-R)J(a/r))^{\rho-1},\\
    \lim_{x\uparrow +\infty}f_v(c(x),J(x))&=\rho\nu \lim 
    _{x\uparrow +\infty}(J'(x))^{-\frac{1-S}{S}} \big((1-R)J(x)\big)^{\frac{(1-S)\rho}{S}}
    ((1-R)J(x))^{\rho-1}.
\end{align*}By (A4) and Proposition \ref{pro-jprime-guji}, we have  
$$(J'(x))^{-\frac{1-S}{S}} \big((1-R)J\big)^{\frac{(1-S)\rho}{S}}
    ((1-R)J(x))^{\rho-1}\sim O\left(x^{\frac{R(1-S)}{S}}x^{\frac{(1-S)\rho(1-R)}{S}} x^{(1-R)(\rho-1)}\right)=O(1),$$
where $\frac{R(1-S)}{S}+\frac{(1-S)\rho(1-R)}{S}+(1-R)(\rho-1)=0$.  Therefore, the drift coefficient in \eqref{eq:fankui-Z} satisfies the locally Lipschitz condition. By Theorem 5.2.5 of \cite{KS91}, equation  \eqref{eq:fankui-Z} admits a unique positive solution $M$. Furthermore, by It\^{o}'s formula, the feedback wealth equation also has a unique solution such that $X_t>\frac{a}{r}$ for all $t\geq0$. Besides, $(\pi,c)\in\mathcal{A}(x)$ follows from the continuity of $X$.   \hfill$\Box$

\begin{remark}\label{rem:underline-z}
 From equation \eqref{eq:fankui-Z}, since $f_v<0$, the comparison theorem for SDEs implies that $M=J'(X)\geq \underline{M}$, where $\underline{M}$ satisfies
\begin{align}\label{eq:undeline-z}
 d\underline{M}_t=(\delta\nu-r)\underline{M}_tdt
    -\frac{\mu-r}{\sigma}\underline{M}_tdB_t, ~~\underline{M}_0=J'(x).   
\end{align}
 This lower bound $\underline{M}$ of the process $M$ admits a lognormal dynamic. This property plays an important role in the final step of the proof of the verification theorem presented below. 
\end{remark}

\subsection*{Proof of Theorem \ref{th:verification}}
The proof is divided into several steps.

{\em Step 1.} We establish that  the value function is the minimal element in the class of all candidate value
functions. That is,
\begin{align}\label{eq: V hat general}
J(x)\geq\sup_{c\in \mathcal{C}(x)} V^c_0.
\end{align}

For a given $x>\frac{a}{r}$,  let $(\pi, c)$ be an arbitrary feasible strategy in
$\mathcal{A}(x)$, and denote $X=(X_{t}^{x,\pi, c})_{t\geq0}$ as the corresponding wealth process. Then, $X_{t}\geq \frac{a}{r}$ for all $t\geq0$. 
Let $Y=(Y_{t})_{t\geq0}$ be the optimal wealth process for Merton's model with the initial wealth of $1$; that is, $Y_{t}=X_{t}^{1,\pi^{(m)}, c^{(m)}}$, where $\pi^{(m)}=\frac{\mu-r}{R\sigma^{2}}$ and $c^{(m)}=\eta Y$. Furthermore, for any $\epsilon>0$, we define the process $Z:=X+\epsilon Y$.  Then $Z_{t}=X_{t}+\epsilon Y_{t}>\frac{a}{r}$ for all $t\geq0$. In particular,  $Z_{0}=x+\epsilon$. The process $Z$ satisfies the following dynamics:
\begin{align*}
dZ_{t}=&d(X_{t}+\epsilon Y_{t})\\
=&\big[rZ_{t}+(X_{t}\pi_{t}+\epsilon \pi^{(m)}Y_{t} )(\mu-r) -(c_{t}+\epsilon \eta Y_{t})\big]dt+(X_{t}\pi_{t}+\epsilon  \pi^{(m)} Y_{t})\sigma dB_{t}\\
=&\big[rZ_{t}+Z_{t}\tilde{\pi}_{t}(\mu-r) -\tilde{c}_{t}\big]dt+Z_{t}\tilde{\pi}_{t}\sigma dB_{t}, 
\end{align*}
where  
\begin{align*}
\tilde{c}_{t}:=c_{t}+\epsilon \eta Y_{t}   \textrm{~~and~~} \tilde{\pi}_{t}:=\frac{X_{t}\pi_{t}+\epsilon  \pi^{(m)} Y_{t} }{Z_{t}}, ~~~\textrm{for all~~} t\geq0.
\end{align*}
It is easy to verify that $(\tilde{\pi}, \tilde{c}_{t})$ is a feasible strategy for initial wealth $x+\epsilon$, i.e., $(\tilde{\pi}, \tilde{c})\in\mathcal{A}(x+\epsilon)$. 

For the given candidate value function $J$, we apply the It\^{o}'s formula to $e^{-\delta\nu t} J\left(Z_t\right)$, $t\geq0$ and obtain
\begin{align*}
d e^{-\delta\nu t} J\left(Z_t\right)
& =e^{-\delta\nu t}\left[-\delta\nu J\left(Z_t\right)+\big(Z_t\left(r+(\mu-r)\tilde{\pi}_t\right)-\tilde{c}_t\big)J^{\prime}(Z_t)+\frac{\sigma^2}{2} \tilde{\pi}_t^2 Z_t^2 J^{\prime\prime}(Z_t)\right] d t\\
&~~~~+e^{-\delta\nu t}Z_t\tilde{\pi}_t\sigma J^{\prime}(Z_t)dB_t\\
&=\Big[e^{-\delta\nu t}L(\tilde{\pi}_t, \tilde{c}_t; Z_t, J(Z_t))-e^{-\delta\nu t}f(\tilde{c}_t,J(Z_t))\Big]dt+dN^{\epsilon}_{t},
\end{align*}  where $N^{\epsilon}_{\cdot}\equiv \int_0^{\cdot} e^{-\delta \nu s} \sigma \tilde{\pi}_s Z_s J^{\prime}(Z_t)dB_s$ is a local martingale and
\begin{align*}
L(\pi, c ; z=x+\epsilon y, v=v(z))=f(c,v)-\delta\nu v+[z(r+(\mu-r)\pi)-c] v_z+\frac{\sigma^2}{2} \pi^2 z^2 v_{z z}.
\end{align*}

By its definition of candidate value function, $J$ satisfies the following HJB equation
\begin{align}\label{eq:Lep}
\sup_{c\geq a, \pi} L(\pi, c ; z=x+\epsilon y, J(z))=0,~~~z>\frac{a}{r}.\end{align}
Maximizing \eqref{eq:Lep} over $\pi \in \mathbb{R}$ and $c \geq a$, we derive the candidate optimal solutions
$$\pi(z):=-\frac{\mu-r}{\sigma^{2}} \frac{J^{\prime}(z)}{z J^{\prime\prime}(z)}  \textrm{~~~and~~~}    c(z):=\max\{a, (J^{\prime}(z))^{-1/S}((1-R)J(z))^{\rho/S}\}. $$
Plugging in $\pi(z)$ and $c(z)$, we get that $L(\pi(z), c(z); z, J(z))=0$ and  
$L\left(\tilde{\pi}_t, \tilde{c}_t ; Z_t, J(Z_t)\right)\leq0$. Then it follows that
\begin{align}\label{eq:Mt-epsilon}
d e^{-\delta\nu t} J\left(Z_t\right)
\leq -e^{-\delta\nu t}f(\tilde{c_t},J(Z_t))dt+dN^{\epsilon}_{t}.
\end{align}
Next, for any two bounded stopping times $\tau_1\leq\tau_2$ and let $$\tau_n=\inf\{s\geq \tau_1: \langle N^{\epsilon}\rangle_s-\langle N^{\epsilon}\rangle_{\tau_1}\geq n\},~~n\geq3.$$ By \eqref{eq:Mt-epsilon}, for all $n\geq 3$, we have 
\begin{align*}
    e^{-\delta\nu \tau_1}J(Z_{\tau_1})
    &\geq e^{-\delta\nu (\tau_2\wedge \tau_n)}J(Z_{\tau_2\wedge\tau_n})+ \int_{\tau_1}^{\tau_2\wedge\tau_n}e^{-\delta\nu s}f(\tilde{c}_s,J(Z_s))ds+N^\epsilon_{\tau_1}-N^\epsilon_{\tau_2\wedge\tau_n}\\
    &=e^{-\delta\nu (\tau_2\wedge \tau_n)}J(Z_{\tau_2\wedge\tau_n})+ \int_{\tau_1}^{\tau_2\wedge\tau_n}e^{-\delta s}f(\tilde{c}_s,e^{-\delta \nu s}J(Z_s))ds+N^\epsilon_{\tau_1}-N^\epsilon_{\tau_2\wedge\tau_n}.
\end{align*}
By taking the conditional expectations on both sides of the last equation and using the optimal sampling theorem, we obtain
\begin{align*} e^{-\delta\nu \tau_1}J(Z_{\tau_1})\geq \mathbb{E}\left[e^{-\delta\nu (\tau_2\wedge \tau_n)}J(Z_{\tau_2\wedge\tau_n})+ \int_{\tau_1}^{\tau_2\wedge\tau_n}e^{-\delta s}f(\tilde{c}_s,e^{-\delta \nu s}J(Z_s))ds~\big|~\mathcal{F}_{\tau_1}\right].
\end{align*}
Since  $J(Z_{\tau_2\wedge\tau_n})=J(X_{\tau_2\wedge\tau_n}+\epsilon Y_{\tau_2\wedge\tau_n})\geq J(\frac{a}{r}),$ then by Fatou's lemma and and the conditional monotone convergence theorem, one derives 
\begin{align*} e^{-\delta\nu \tau_1}J(Z_{\tau_1})\geq \mathbb{E}\left[e^{-\delta\nu \tau_2}J(Z_{\tau_2})+ \int_{\tau_1}^{\tau_2}e^{-\delta s}f(\tilde{c}_s,e^{-\delta \nu s}J(Z_s))ds~\big|~\mathcal{F}_{\tau_1}\right].
\end{align*}
Moreover, 
\begin{align*}
\liminf_{t\rightarrow+\infty}\mathbb{E}\left[e^{-\delta\nu t}J(Z_t)\right]\geq \liminf_{t\rightarrow+\infty}\mathbb{E}\left[e^{-\delta\nu t}J(X_t)\right]\geq \lim_{t\rightarrow+\infty} e^{-\delta\nu t}J(\frac{a}{r})=0.
\end{align*}
Therefore, the utility process 
$e^{-\delta\nu t}J(Z_{t})$, $t\geq0$, is a supersolution in the sense of \cite{HHJ23b} (Definition 5.3) and \cite{SS16} for the $(f,\tilde{c})$, where $\tilde{c}=c+\eta\epsilon Y$.

We now show that, for the candidate value function $J$, 
\begin{align}\label{eq:Jhshagnjie}
   J(x+\epsilon)\geq \sup_{c\in \mathcal{C}(x)} V^c_0,~~~x>\frac{a}{r}. 
\end{align}
For $R<1$, since $f(c,v)$ is increasing with respect to $c$. Therefore, $e^{-\delta\nu t}J(Z_{t})$, $t\geq0$,  is also a supersolution for  $(f,c)$.
By \cite{HHJ23b}, Theorem 6.5, which asserts that $V^c$ is the minimal supersolution to $(f,c)$, we conclude that $J(x+\epsilon )\geq V^c_0$ for all $c\in\mathcal{C}(x)$. 

For $R>1$, since $\tilde{c}^{1-S}\leq (\eta \epsilon)^{1-S} Y^{1-S}$, then there exists a utility process $V^{\tilde{c}}$ such that uniformly integrable. Then, by Proposition 6.8 and Theorem 6.5 in \cite{HHJ23b}, it follows that $J(x+\epsilon)\geq V_0^{\tilde{c}}\geq V_0^{c}$ for all $c\in\mathcal{C}(x)$, and \eqref{eq:Jhshagnjie} holds.

Finally, setting $\epsilon \downarrow 0$, we obtain that  \eqref{eq: V hat general} holds for all $x>\frac{a}{r}$. 

{\em Step 2.} For the candidate value function $J$, let $(\pi^*, c^*)$ be the feedback control induced by \eqref{eq:h-strategy}, and let $X^{*}$ be the corresponding feedback wealth process. We prove the verification result under a so-called strong transversality condition \begin{equation}
\label{eq:strongtransversality}
\lim_{t\rightarrow+\infty}\mathbb{E}\left[e^{-\delta\nu t} J(X_{t}^{*})\right]=0. 
\end{equation}

For any $x>\frac{a}{r}$, for the candidate value function $J$, it follows from Proposition \ref{pro: SDE} that $X^*_t>\frac{a}{r}$. 
Similar to the proof procedure of step 1, 
applying the It\^{o}'s formula to $e^{-\delta\nu t} J\left(X^*_t\right)$, $t\geq0$, we obtain
\begin{align*}
d e^{-\delta\nu t} J\left(X^*_t\right)&=\Big[e^{-\delta\nu t}L(\pi^*(X^*_t), c^*(X^*_t); X^*_t, J(X^*_t))-e^{-\delta\nu t}f(c^*(X_t^*),J(X^*_t))\Big]dt+dN_{t}\nonumber\\
&=-e^{-\delta\nu t}f(c^*(X_t^*),J(X^*_t))dt+dN_{t},
\end{align*}  where $N_{\cdot}\equiv \int_0^{\cdot} e^{-\delta \nu s} \sigma \pi^*(X^*_s) X^*_s J^{\prime}(X^*_t)dB_s$. 
Choosing a suitable stopping time $\tau$ and taking the expectation on the both sides, it implies that 
\begin{align*} J(x)= \mathbb{E}\left[e^{-\delta\nu \tau}J(X^*_{\tau})+ \int_{0}^{\tau}e^{-\delta s}f(c^*(X_s^*),e^{-\delta \nu s}J(X^*_s))ds\right].
\end{align*}
Letting $\tau \rightarrow \infty$, we obtain 
$$J(x)=V_0^{c^*} =\mathbb{E}\left[\int_{0}^{+\infty}e^{-\delta s}f(c^*(X_{s}^{*}), V_s^{c^*})ds\right],$$
if the above strong transversality condition \eqref{eq:strongtransversality} holds. Combining the result of Step 1, we conclude that the strategy $(\pi^*, c^*)$ attains the value function.


{\em Step 3.} We prove that the strong transversality condition \eqref{eq:strongtransversality} holds, thereby completing the proof.

When $R>1$, by (A1), $J$ is negative. Moreover, 
\begin{align*}
0\geq\limsup_{t\rightarrow+\infty}\mathbb{E}\left[e^{-\delta\nu t} J(X_{t}^{*})\right]\geq\liminf_{t\rightarrow+\infty}\mathbb{E}\left[e^{-\delta\nu t} J(X_{t}^{*})\right]\geq \liminf_{t\rightarrow+\infty}\mathbb{E}\left[e^{-\delta\nu t} J(\frac{a}{r})\right]=0.
\end{align*}Thus, the strong transversality condition holds.  

We next consider the case $0 < R < 1$. In this situation, by (A4), Proposition \ref{pro-jprime-guji} and  Remark \ref{rem:underline-z}, we have
\begin{align*}
  0\leq \limsup_{t\rightarrow\infty}\mathbb{E}\left[e^{-\delta\nu t} J(X^*_{t})\right]&\leq \limsup_{t\rightarrow\infty}\frac{\eta^{-\nu S}}{1-R}
  \mathbb{E}\left[e^{-\delta\nu t} (X_t^*)^{1-R}\right]\\
  &= \limsup_{t\rightarrow\infty}\frac{\eta^{-\nu S}}{1-R}
  \mathbb{E}\left[I_{X_t^*\leq x_0}e^{-\delta\nu t} (X_t^*)^{1-R}+I_{X_t^*> x_0}e^{-\delta\nu t} (X_t^*)^{1-R}\right]\\
  &\leq \limsup_{t\rightarrow\infty}\frac{\eta^{-\nu S}}{1-R}x_0^{1-R}e^{-\delta\nu t}+ k_2^{-1+1/R}\frac{\eta^{-\nu S}}{1-R}\mathbb{E}\left[e^{-\delta\nu t}\Big(J'(X^*_t)\Big)^{1-\frac{1}{R}}\right]\\
  &\leq \limsup_{t\rightarrow\infty}k_2^{-1+1/R}\frac{\eta^{-\nu S}}{1-R}\mathbb{E}\left[e^{-\delta\nu t}(\underline{M}_t)^{1-\frac{1}{R}}\right].
  \end{align*}

By virtue of \eqref{eq:undeline-z}, we obtain
\begin{align*}
  \mathbb{E}\left[e^{-\delta\nu t}(\underline{M}_t)^{1-\frac{1}{R}}\right]
  &\leq (\underline{M}_0)^{1-1/R} e^{-\delta\nu t}\mathbb{E}\left[e^{(\delta \nu-r-\kappa)(1-\frac{1}{R})t-\frac{\mu-r}{\sigma}(1-\frac{1}{R})B_t}\right]\\
  &=(\underline{M}_0)^{1-1/R} e^{-\delta\nu t} e^{(\delta \nu-r-\kappa)(1-\frac{1}{R})t+\kappa(1-\frac{1}{R})^2t}=(\underline{M}_0)^{1-1/R}e^{-\frac{\eta \nu S}{R} t}.
\end{align*}This completes the proof of the strong transversality condition for $\eta>0$. \hfill$\Box$


\subsection*{Proof of Theorem \ref{th:two-region}}
We first observe that $(1-R)J > 0$, and by
equation \eqref{eq:bound}, $((1-R)J)^{\frac{\rho}{S}}
$ remains bounded when $x \downarrow \frac{a}{r}$. Furthermore, by Proposition \ref{pro:J-daoshu}, 
$\lim_{x \downarrow \frac{a}{r}}J'(x) = +\infty$. 
Consequently, by its definition, $\mathcal{B} =\big\{ ((1-R)J)^{\frac{\rho}{S}} (J')^{-\frac{1}{S}}<a\big\}$, there exists a positive number $b > \frac{a}{r}$ such that $(\frac{a}{r}, b) \subseteq {\cal B}$. Moreover, $b\neq+\infty$; otherwise, the constraint region $\mathcal{B}$ should be $(\frac{a}{r}, +\infty)$, implying that the consumption rate always takes the minimal rate $a$. In this case,  $J(x)\equiv \delta^{-\nu}\frac{a^{1-R}}{1-R}$, $\forall x\geq {a\over r}$, which contradicts Proposition \ref{prop:strictly increasing} that $J(x)$ is strictly increasing.


Since $c^{*}=\max\Big\{a, (J'(x))^{-1/S} \big((1-R)J(x)\big)^{\frac{\rho}{S}}\Big\}$,  it suffices to show that $(J'(x))^{-1} \big((1-R)J(x)\big)^{\rho}$ increases with respect to $x$ on $(\frac{a}{r},+\infty)$. A  straightforward calculation yields 
\begin{equation*}
((J'(x))^{-1} \big((1-R)J(x)\big)^{\rho})'={((1-R)J)^{\rho-1}\over (J'(x))^2}\cdot (1-R)\Big(\rho(J'(x))^2-J(x)\cdot J''(x)\Big).
\end{equation*}
Therefore, it remains to prove that $(1-R)[\rho(J'(x))^2-J(x)\cdot J''(x)]>0$. 

We analyze the problem by considering the following two cases separately.

{\em Case 1. $R > 1$.}  In this case, since $J(x)<0$, combined with $\rho<0$ and $J''(x)<0$, it follows directly that $(1-R)[\rho(J'(x))^2-J(x)\cdot J''(x)]>0$.

{\em Case 2. $R < 1$.} In this situation, since $J(x)>0$, $\forall x>\frac{a}{r}$, it is equivalent to show that 
\begin{align}
 \label{eq:region}
 {(J'(x))^2\over J(x)\cdot J''(x)}>{1\over \rho},~~~x>\frac{a}{r}.
\end{align}
Since $J$ is $C^2$, it suffices to verify (\ref{eq:region}) separately on $\mathcal{B}$ and $\mathcal{U}$.

In the constraint region $\mathcal{B}$, by (\ref{eq:J in B}), we have
\begin{align*}
\kappa {J'^2\over J''}&={a^{1-S}\over 1-S}((1-R)J)^{\rho}-aJ'+rxJ'-\delta\nu J\\
\Longrightarrow {J'^2\over J J''}&={1\over \kappa}\big[{a^{1-S}\over 1-S}{((1-R)J)^{\rho}\over J}+(rx-a){J'\over J}-\delta\nu \big]>-{\delta \nu \over \kappa}.
\end{align*}

Similarly, in the unconstrainted region $\mathcal{U}$, by (\ref{eq:J in U}), we have
\begin{align*}
\kappa {J'^2\over J''}&={S\over 1-S}((1-R)J)^{\rho \over S}(J')^{1-{1\over S}}+rxJ'-\delta\nu J\\
\Longrightarrow {J'^2\over J J''}&={1\over \kappa}\big[{S\over 1-S}{((1-R)J)^{\rho \over S}(J')^{1-{1\over S}}\over J}+rx{J'\over J}-\delta\nu  \big]>-{\delta \nu \over \kappa}.
\end{align*}
Therefore, (\ref{eq:region}) follows from the condition $\kappa+\rho\delta\nu>0$. The proof is complete.
\hfill$\Box$

\end{document}